\documentclass[a4paper]{amsart}
\usepackage{amsmath,amssymb}
\usepackage[dvipsnames]{xcolor}
\usepackage{xcolor}
  \def\red#1{\textcolor{red}{#1}} 
  
\newtheorem{theorem}{Theorem}[section]
\newtheorem{proposition}[theorem]{Proposition}
\newtheorem{corollary}[theorem]{Corollary}

\newtheorem{algorithm}[theorem]{Algorithm}

\newtheorem{preremark}[theorem]{Remark}
\newtheorem{predefinition}[theorem]{Definition}
\newtheorem{preexample}[theorem]{Example}
\newtheorem{prenotation}[theorem]{Notation}
\newtheorem{preconjecture}[theorem]{Conjecture}
\newtheorem{preassumption}[theorem]{Assumption}

\newenvironment{remark}{\begin{preremark}\rm}{\end{preremark}}
\newenvironment{definition}{\begin{predefinition}\rm}
{\end{predefinition}}
\newenvironment{example}{\begin{preexample}\rm}{\end{preexample}}

\allowdisplaybreaks
\def\ZZ{\mathbb{Z}}
\def\QQ{\mathbb{Q}}

\newcommand{\M}{{\mathfrak{M}}}

\newcommand{\m}{{\mathfrak{m}}}

\let\epsilon=\varepsilon

\def\phi{{\varphi}}
\let\Psi=\varPsi
\let\Phi=\varPhi
\let\theta=\vartheta
\let\rho=\varrho

\def\LT{\mathop{\rm LT}\nolimits}

\def\GL{\mathop{\rm GL}\nolimits}

\def\Mat{\mathop{\rm Mat}\nolimits}

\def\Supp{\mathop{\rm Supp}\nolimits}
\def\Spec{\mathop{\rm Spec}\nolimits}
\def\Proj{\mathop{\rm Proj}\nolimits}

\def\Cot{\mathop{\rm Cot}\nolimits}

\def\Syz{\mathop{\rm Syz}\nolimits}
\def\Ker{\mathop{\rm Ker}\nolimits}

\def\lin{{\mathop{\rm lin}\nolimits}}

\def\cvec{{\mathop{\rm cvec}\nolimits}}
\def\cmat{{\mathop{\rm cmat}\nolimits}}

\newcommand{\Lin}{\mathop{\rm Lin}\nolimits}

\newcommand{\Xz}{X^\circ}
\newcommand{\Xp}{X^\plus}

\let\To=\longrightarrow

\def\TTTo#1\mathop{\xrightarrow{\hspace*{1cm}}^{^{\mkern-70mu (#1}}}

\def\tfrac #1#2{{\textstyle\frac{#1}{#2}}}
\def\tsum_#1^#2{{\textstyle\sum\limits_{#1}^{#2}}}

\def\plus{{\scriptscriptstyle +}}

\definecolor{red}{rgb}{1.0, 0.0, 0.0}

\def\cocoa{\mbox{\rm
  C\kern-.13em o\kern-.07 em C\kern-.13em o\kern-.15em A}}
\def\apcocoa{\mbox{\rm
A\kern-0.13em p\kern -0.07em C\kern-.13em o\kern-.07 em C\kern-.13em
o\kern-.15em A}}

\begin{document}

\title{Elimination by Substitution}

\author{Martin Kreuzer}
\address{Fakult\"at f\"ur Informatik und Mathematik, Universit\"at
Passau, D-94030 Passau, Germany}
\email{martin.kreuzer@uni-passau.de}

\author{Lorenzo Robbiano}
\address{Dipartimento di Matematica, Universit\`a di Genova,
Via Dodecaneso 35,
I-16146 Genova, Italy}
\email{lorobbiano@gmail.com}

\date{\today}

\begin{abstract}
Let $K$ be a field and $P=K[x_1,\dots,x_n]$. 
The technique of elimination by substitution is based on discovering a coherently
$Z=(z_1,\dots,z_s)$-separating tuple of polynomials $(f_1,\dots,f_s)$ in an ideal~$I$, i.e., 
on finding polynomials such that $f_i = z_i - h_i$ with $h_i \in K[X{\setminus} Z]$. Here we
elaborate on this technique in the case when~$P$ is non-negatively graded. The existence
of a coherently $Z$-separating tuple is reduced to solving several $P_0$-module membership problems.
Best separable re-embeddings, i.e., isomorphisms $P/I \longrightarrow K[X{\setminus Z}] /
(I \cap K[X {\setminus} Z])$ with maximal $\#Z$, are found degree-by-degree. 
They turn out to yield optimal re-embeddings in the
positively graded case. Viewing $P_0 \longrightarrow P/I$ as a fibration over an affine space, we
show that its fibers allow optimal $Z$-separating re-embeddings, and we provide a criterion for
a fiber to be isomorphic to an affine space. In the last section we introduce a new technique based on
the solution of a unimodular matrix problem which enables us to construct automorphisms of~$P$ 
such that additional $Z$-separating re-embeddings are possible. 
One of the main outcomes is an algorithm which allows us to explicitly compute a homogeneous
isomorphism between $P/I$ and a non-negatively graded polynomial ring if $P/I$ is regular.
\end{abstract}



\maketitle


%
%

\section{Introduction}
\label{sec1}

The computation of elimination ideals is an important task in computer algebra,
as it corresponds geometrically to calculating the image of a projection.
Classically, elimination was performed by computing {\it resultants}. 
This method typically requires the calculation of the determinants of huge polynomial matrices and is
impracticable when the number of indeterminates is large. In the last 50 years,
the preferred method for performing elimination has been the computation of 
{\it Gr\"obner bases} with respect to elimination term orderings. Despite their high
worst-case complexity, algorithms for computing Gr\"obner bases work well in many cases
if the number of indeterminates is not too large. One disadvantage of elimination Gr\"obner
basis calculations is that the polynomials generating the elimination ideal are found near
the end of the computation. Moreover, they can be rather numerous and big, because one calculates
a Gr\"obner basis of the elimination ideal, not just a system of generators.

In the last years, we have been working with the ideals defining border basis schemes 
(e.g., see~\cite{KR3}, \cite{KLR0}). They are explicitly given and allow us to study an 
important class of schemes, but they involve a large number of indeterminates. It is therefore
desirable to embed them into lower-dimensional polynomial rings, 
i.e., to eliminate as many indeterminates as possible.
Unfortunately, the usual elimination methods break down in all but the easiest examples. 
Therefore, in a series of papers (see~\cite{KLR1}, \cite{KLR2}, \cite{KLR3}) 
together with Le N.\ Long, we developed the method of {\it $Z$-separating re-embeddings}
which is based on a way to compute {\it elimination by substitution}. Its advantages are that
substitution can be calculated rather efficiently and that the result is a system of generators of the
elimination ideal which is simpler (and possibly much simpler) than a Gr\"obner basis.
Its drawbacks are that it does not apply in general, and that it
may eliminate only some of the indeterminates which could be eliminated using other methods
(for instance, see Examples~\ref{ex-twodeg} and~\ref{ex-similCrachiolaCont}).
Notice that, in the final section, we introduce such additional methods, as well.

In this paper we develop this technique of elimination by substitution further in the case
of homogeneous ideals in non-negatively graded polynomial rings. This setting was also suggested
to us by ideals defining so-called maxdeg border basis schemes. However, it has turned out to
be very rich and interesting in its own right. Geometrically, the residue class ring of a non-negatively
graded polynomial ring modulo a homogeneous ideal which is generated in positive degrees corresponds to
a fibration over an affine space whose fibers are affine spaces under mild assumptions. 
Not only were we able to construct new methods for finding
best possible, and sometimes optimal $Z$-separating re-embeddings of homogeneous ideals, but we 
also constructed a class of isomorphisms of polynomial rings which is based on the unimodular matrix problem
and which allows us to eliminate further indeterminates by substitution. Moreover, we found
criteria under which we can effectively show that the above fibration is globally isomorphic to an affine space.
\smallskip

Now we describe the contents of this paper step-by-step. In Section~2 we first recall the notions
of $Z$-separating and coherently $Z$-separating tuples of polynomials and then formulate the basic
technique for elimination by substitution. More precisely, given a polynomial
ring $P=K[x_1,\dots,x_n]$ over a field~$K$ and a tuple of distinct indeterminates $Z=(z_1,\dots,z_s)$
contained in $X=(x_1,\dots,x_n)$ we say that a tuple $F=(f_1,\dots,f_s)$ is {\it coherently $Z$-separating} 
if $f_i = z_i - h_i$ with $h_i \in K[X{\setminus} Z]$ for $i=1,\dots,s$. When we find such a tuple in an ideal~$I$, 
we can compute the elimination ideal $I\cap K[X{\setminus} Z]$ by substitution,
i.e., we can compute a {\it $Z$-separating re-embedding} $\phi:\; P/I \longrightarrow
K[X{\setminus}Z] / (I \cap K[X{\setminus}Z]$ (see Proposition~\ref{prop-ElimViaSubst})
without resorting to the computation of a Gr\"obner basis.
If $I$ is contained in~$\M=\langle x_1,\dots,x_n\rangle$, candidates for suitable tuples~$Z$ can be recognized
using the minors of the coefficient matrix of the linear parts of a system of generators of~$I$ 
(see Proposition~\ref{prop-linindep}), but there remains the major problem of finding out whether
a coherently $Z$-separating tuple of polynomials exists in~$I$.

In order to go further, we then introduce the main setting of this paper in Section~3.
We assume that~$P$ is non-negatively graded by $W=(w_1,\dots,w_m)$ and that~$I$
is a $W\!$-homogeneous ideal in~$P$ such that $I\cap P_0 = \{0\}$.
Let us point out right away that most of the results become substantially
easier and are partially known if the grading by~$W$ is positive. So, our main focus is on the case
when $P_0$ is a polynomial ring in one or more indeterminates. Then the geometric interpretation
of the $P_0$-algebra $P/I$ is that it yields a fibration over an affine space. 

In this setting, a first preprocessing step is Algorithm~\ref{alg-FindSepPoly} which produces the list of all 
indeterminates~$z_i$ such that~$I$ contains a $z_i$-separating polynomial.
Then we note that every polynomial $g\in \langle X\rangle$ has a representation
$g = p_1 z_1 + \cdots + p_s z_s + h$ with $p_i \in P_0$ and such that no term $p\,z_i$ with 
$p\in P_0$ is in $\Supp(h)$. 
Here the tuple $\cvec_Z(g)=(p_1,\dots,p_s)$ is called the {\it $Z$-coefficient vector}
of~$g$. Our first main result is Theorem~\ref{thm-unitvec} which, given a tuple~$Z$, 
allows us to check whether~$I$ contains a coherently $Z$-separating tuple. If it does, then Algorithm~\ref{alg-findF} 
helps us to find such a tuple. The key point is that everything reduces to solving several $P_0$-submodule
membership problems for the $P_0$-module generated by the $Z$-coefficient vectors of a homogeneous system
of generators of~$I$. 

In Section~4 we then apply this technique to find {\it best separable re-embeddings}
in the non-negatively graded case. These are separating re-embeddings for which~$\#Z$
is as large as possible. A best separable re-embedding is not necessarily {\it optimal}, i.e., it does not
necessarily reach the minimal embedding dimension, but it is the best we can
do using the elimination-by-substitution technique, and in many cases it is already quite good.
We start by considering the case of tuples~$Z$, all of whose entries have the same $W\!$-degree~$d$.
Algorithm~\ref{alg-BestZSepInDegd} computes a best separable tuple~$Z$ in this case, and 
Remark~\ref{rem-allbest}
provides a variant which finds {\it all} best separating tuples~$Z$.
For the general case, we show in Proposition~\ref{prop-deg-by-deg} that homogeneous separating
tuples can be computed degree-by-degree, and in Algorithm~\ref{alg-BestPossible} we combine everything
to find best separating tuples for homogeneous ideals in the non-negatively graded case.
When we specialize this to the positively graded case in Corollary~\ref{cor-optimalgraded},
we find that our algorithms automatically produce optimal re-embeddings, recovering a version
of a result of M.\ Roggero and L.\ Terracini (cf.~\cite{RT}).

Coming back to the point of view that $P_0 \longrightarrow P/I$ corresponds to  a fibration over an affine space,
we study the fibers of this map in Sections~5. For this purpose we adapt our notation slightly:
now the indeterminates of degree zero are called $a_1,\dots,a_m$, and the indeterminates of positive degree
are the entries of $\Xp=(x_1,\dots,x_{n-m})$. Thus the base space for the fibration is the affine space
defined by $P_0 = K[a_1,\dots,a_m]$. For a point $\Gamma = (c_1,\dots,c_m) \in K^m$, we let $\M_\Gamma =
\langle a_1-c_1,\dots,a_m-c_m, x_1,\dots,x_{n-m}\rangle$ be the corresponding maximal ideal of~$P$.
We call the ideal $I_\Gamma \subseteq K[\Xp]$ obtained from~$I$ via the substitutions $a_i \mapsto c_i$
the {\it fiber ideal} over~$\Gamma$. After noting that coherently separating tuples remain coherently 
separating in the fibers (see Proposition~\ref{prop-CohSepInFibers}), we prove in Proposition~\ref{prop-OptFibers}
that there {\it always} exist tuples of indeterminates~$Z$ which actually yield optimal re-embeddings 
of the fibers. This is one of the main perks of working in the non-negatively graded setting. 
Our second main result in this section is Theorem~\ref{thm-freefiber} which says (among others)
that the regularity of the local ring $(P/I)_{\M_\Gamma/I}$ implies that $K[\Xp]/I_\Gamma$ is isomorphic
to a polynomial ring, i.e., that the fiber over~$\Gamma$ is isomorphic to an affine space
of dimension $\dim((P/I)_{\M_\Gamma/I})-m$. Moreover, if~$P/I$ is a regular ring, then it is an integral
domain and the fibers have dimension $\dim(P/I)-m$ (see Proposition~\ref{prop-isprime}).
In the language of Algebraic Geometry, this says that the morphism $\Spec(P/I) \longrightarrow \Spec(P_0)$
is an $\mathbb{A}^k$-fibration, where $k=\dim(P/I)-\dim(P_0)$.

In Section~6 we move from studying the individual fibers of the morphism $\Spec(P/I) \longrightarrow \Spec(P_0)$
to looking at the entire family. The question whether an $\mathbb{A}^k$-fibration is isomorphic to an affine space
is, in general, a hard problem (see~\cite{Gup}, Section~4, Question~3).
Here we bring in a technique which allows us to find $K$-algebra
automorphisms $\phi:\; P \longrightarrow P$ such that~$\phi(I)$ has a $Z$-separating re-embedding, even if~$I$
doesn't. This technique is based on the {\it unimodular matrix problem (UMP)} which says that, given
a matrix $A\in \Mat_{k,n-m}(P_0)$ whose maximal minors generate the unit ideal, we can compute a
matrix $B\in \GL_{n-m}(P_0)$ such that $A\cdot B = (I_k \mid 0)$ where $I_k$ denotes the identity matrix of
size~$k$. The effective solution of this task corresponds to {\it unimodular matrix completion}, resp.\ 
the effective version of the Quillen-Suslin theorem (cf.~\cite{Kun}, IV.3.15), and has been implemented in 
several computer algebra systems. 

Based on a suitably adapted version of this method (see Algorithm~\ref{alg-completion}),
we provide in Algorithm~\ref{alg-UMC-reemb} a general method for re-embedding the ideal~$I$ if
it contains homogeneous polynomials $g_i = h_{i1} x_1 + \cdots + h_{i\, n-m} x_{n-m} + \tilde{h}_i$
with $h_{ij}\in P_0$ and $\tilde{h}_i \in \langle \Xp\rangle^2$ such that the maximal minors of the matrix of coefficients
$A = (h_{ij}) \in \Mat_{k.nm}(P_0)$ generate the unit ideal. More precisely, we first construct
an automorphism of~$P$ using the UMP solution technique and then perform a $Z$-separating re-embedding.

Finally, we apply this algorithm to the case when $P/I$ is a regular ring. In this case the hypotheses
of Algorithm~\ref{alg-UMC-reemb} are satisfied and we obtain a
homogeneous \hbox{$K$-algebra} isomorphism between $P/I$ and a non-negatively graded polynomial ring~$\widehat{P}$
(see Theorem~\ref{thm-Free}). In conclusion, Algorithm~\ref{alg-free} lays out the steps to compute this isomorphism
$P/I \cong \widehat{P}$ effectively.

Overal, let us emphasize that the re-embeddings found in Sections~3,~4, and~5 could, in principle, be calculated 
using Gr\"obner basis methods. However, as mentioned before, this approach will fail in all but the simplest cases. 
On the other hand, the re-embeddings found in Section~6 are completely outside the
scope of Gr\"obner basis techniques and require, in general, the effective solution of the unimodular matrix problem.

For the notations and basic definitions used in this paper, we refer to our books~\cite{KR1} and~\cite{KR2}.
All examples have been calculated using the computer algebra system \cocoa\ (see~\cite{CoCoA})
and their source files are available from the authors upon request.

\bigskip\bigbreak
%
%

\section{Coherently Separating Tuples and Elimination by Substitution}
\label{sec2}

In this paper we let $K$ be a field, let $P=K[x_1,\dots,x_n]$ be a polynomial ring over~$K$, 
let $\M=\langle x_1,\dots,x_n \rangle$ be the maximal ideal generated by the indeterminates, and let
$X=(x_1,\dots,x_n)$ be the tuple of indeterminates in~$P$.

In  the first part of the paper we use the method of $Z$-separating re-embeddings, as developed 
in~\cite{KLR1}, \cite{KLR2}, and~\cite{KLR3}. Let us recall the central definitions.

\begin{definition}{\bf (Separating and Coherently Separating Tuples)}\label{def-sepindets}\,\\
Let $Z=(z_1,\dots,z_s)$ be a tuple of distinct indeterminates in $X=(x_1,\dots,x_n)$.
\begin{enumerate}
\item[(a)] Let $i\in \{1,\dots, s\}$. A polynomial $f\in P$ is called {\bf $z_i$-separating}
if it is of the form $f=z_i - h$ with a polynomial~$h$ such that $z_i$ divides
no term in $\Supp(h)$.

\item[(b)] A tuple of polynomials $F=(f_1,\dots,f_s) \in P^s$ is called
{\bf $Z$-separating} if there exists a term ordering~$\sigma$ such that
$\LT_\sigma(f_i)=z_i$ for $i=1,\dots, s$.

\item[(c)] A tuple of polynomials $F=(f_1,\dots,f_s) \in P^s$ is called
{\bf coherently $Z$-separating} if the polynomials~$f_i$ are of the form $f_i = z_i - h_i$
with polynomials~$h_i$ such that $z_i$ divides no term in $\Supp(h_j)$
for $i,j \in \{1,\dots,s\}$.

\item[(d)] Let $I$ be an ideal in~$P$. The tuple~$Z$ is called a  {\bf separating tuple of indeterminates} 
for~$I$ if~$I$ contains a $Z$-separating tuple of polynomials.

\end{enumerate}
\end{definition}

Next we recap some properties of separating tuples (for further details, see~\cite{KLR2}).

\begin{remark}\label{rem-sepandcoherentlysep}
Let $Z$ be a tuple of distinct indeterminates in~$X$.
\begin{enumerate}
\item[(a)] Given a $Z$-separating tuple $F=(f_1,\dots,f_s)$, the polynomial~$f_i$
is $z_i$-separating for $i=1,\dots,s$.

\item[(b)] A coherently $Z$-separating tuple is $Z$-separating.

\item[(c)] A $Z$-separating tuple~$F$ is a minimal monic 
$\sigma$-Gr\"obner basis of $\langle F \rangle$.

\item[(d)] A coherently $Z$-separating tuple~$F$ is the reduced 
$\sigma$-Gr\"obner basis of $\langle F\rangle$.

\item[(e)] Given a $Z$-separating tuple~$F$, a coherently $Z$-separating
tuple is obtained by interreducing~$F$. For our setting, explicit instructions are given
in Algorithm~\ref{alg-findCohF}.

\end{enumerate}
\end{remark}

Given a coherently $Z$-separating tuple of polynomials in an ideal~$I$, we can rewrite other 
polynomials in~$P$ as follows.

\begin{definition}
Let $Z=(z_1,\dots,z_s)$ be a tuple of distinct indeterminates in~$X$, and let
$F=(f_1,\dots,f_s)$ be a coherently $Z$-separating tuple of polynomials in an ideal~$I$,
where $f_i = z_i - h_i$  for $i=1,\dots,s$.
\begin{enumerate}
\item[(a)] For every polynomial $g\in P$ and every $i\in \{1,\dots,s\}$,
replace each occurrence of~$z_i$ in~$g$ by~$h_i$. Then the resulting polynomial 
$\hat{g}$ is said to be obtained by {\bf rewriting~$g$ using~$F$}.

\item[(b)] Given a tuple of polynomials $G=(g_1,\dots,g_r)$ in~$P$, we let
$\widehat{G} = (\hat{g}_1,\dots,\hat{g}_r)$ and call it the tuple obtained by
{\bf rewriting~$G$ using~$F$}.

\end{enumerate}
\end{definition}

The power of coherently $Z$-separating tuples of polynomials derives from the fact that
they allow us to perform elimination of~$Z$ as follows.

\begin{proposition}{\bf (Computing Elimination By Substitution)}
\label{prop-ElimViaSubst}\\
Let $K$ be a field, let $P=K[x_1,\dots,x_n]$, let $I=\langle G\rangle$ be an ideal in~$P$
generated by a tuple of polynomials $G=\langle g_1,\dots,g_r)$, let $Z=(z_1,\dots,z_s)$
be a tuple of distinct indeterminates in~$X$, and let $F=(f_1,\dots,f_s)$ be a coherently
$Z$-separating tuple of polynomials in~$I$.

Then the tuple $\widehat{G}$ obtained by rewriting~$G$ using~$F$ is a system of generators
of the elimination ideal $I\cap K[X\setminus Z]$.
\end{proposition}

\begin{proof}
For $i=1,\dots, s$, let $f_i= z_i - h_i$ with $h_i\in P$.
Since $z_i-h_i \in I$ for $i=1,\dots,s$, the result  of applying the substitution
$z_i\mapsto h_i$ to an element $g_j \in I$ yields an element in~$I$. 
Hence the results of applying all these substitutions to the elements 
in~$G$ yields elements in $I\cap K[X\setminus Z]$.
Clearly, every element of $I \cap K[X\setminus Z]$ is obtained in this way, 
for instance from itself, and the proof is complete.
\end{proof}

The following example shows this method at work. Since it is quite small,
the result can be verified using a Gr\"obner basis computation.

\begin{example}\label{ex-ElimBySubst}
Let $P = \QQ[a, w, x, y]$, and let $I =\langle g_1, g_2, g_3 \rangle$, where $g_1 = 3x -ay$, 
$g_2 = y -a^3y$, and $g_3 = ax - a^2w^2$. Now we consider the polynomials
\begin{align*}
f_1 &= (-\tfrac{1}{3}a^3 +\tfrac{1}{3}) g_1 +\tfrac{1}{3}a g_2 +a^2 g_3 \;=\;  x - a^4w^2\in I\\
f_2 &= -a^2 g_1 + g_2 + 3a g_3  \;=\; y - 3a^3w^2 \in I
\end{align*}
and notice that $F=(f_1,f_2)$ is coherently $Z$-separating for $Z=(x,y)$.
Later we will see how to compute tuples like $(f_1,f_2)$ (cf.\ Algorithm~\ref{alg-findF}).

The substitutions $x \mapsto a^4 w^2$ and $y \mapsto 3a^3 w^2$ in $G=(g_1,g_2,g_3)$
yield the triple $\widehat{G} = (0,\; -3 a^6 w^2 + 3a^3 w^2,\; a^5 w^2 - a^2w^2)$ 
which generates the ideal $I \cap \QQ[a,w] = \langle a^5 w^2 - a^2 w^2 \rangle$.
\end{example}

Geometrically, elimination corresponds to projection. 
Thus $Z$-separating tuples allow us to re-embed ideals in the following sense.

\begin{definition}
Let $I$ be an ideal in~$P$, let $Z=(z_1,\dots,z_s)$ be a tuple of distinct indeterminates in~$X$, 
and let $F=(f_1,\dots,f_s)$ be a coherently $Z$-separating tuple in~$I$.
For $i=1,\dots,s$, write $f_i = z_i - h_i$ with $h_i\in P$.
Then the $K$-algebra homomorphism
$$
\Phi:\; P/I \longrightarrow K[X\setminus Z] / (I\cap K[X\setminus Z]) 
$$
defined by $z_i \mapsto h_i$ for $i=1,\dots, s$ and $x_j \mapsto x_j$ for
$x_j \notin Z$ is an isomorphism and is called a {\bf $Z$-separating re-embedding} of~$I$.
\end{definition}

Recall that $\Lin_\M(I)$ denotes the {\bf linear part} of~$I$
(cf.~\cite{KLR1}, Section~1). Here {\it linear}\/ refers to the standard degree,
and $\Lin_\M(I)$ is the $K$-vector space generated by the linear parts
of the polynomials in~$I$. Moreover, if $L=(\ell_1, \dots, \ell_r)$ is  a tuple of
linear forms in~$P$ and if we write $\ell_i = a_{i1} x_1 + \cdots + a_{in} x_n$
with $a_{ij} \in K$ for $i=1,\dots, r$ and $j=1,\dots,n$, then the matrix 
$A=(a_{ij}) \in \Mat_{r,n}(K)$ is called the {\bf coefficient matrix} of~$L$.

Using the linear part of~$I$, candidates for separating tuples~$Z$ for~$I$ can be recognized as follows.

\begin{proposition}\label{prop-linindep}
Let $g_1,\dots,g_r\in \M$, let $I=\langle g_1,\dots,g_r\rangle$, let
$Z=(z_1,\dots,z_s)$ be a tuple of distinct indeterminates, and let
let $A \in \Mat_{r,n}(K)$ be the coefficient matrix of $(\Lin_\M(g_1),\dots, \Lin_\M(g_r))$.

If $Z$ is a separating tuple for~$I$, then the submatrix~$A'\in\Mat_{r,s}(K)$ of~$A$ 
which consists of the columns corresponding to the indeterminates in~$Z$ has rank~$s$. 
In particular, we have $s \le \dim_K(\Lin_\M(I))$.
\end{proposition}

\begin{proof}
Let $F=(f_1,\dots,f_s)$ be a $Z$-separating tuple in~$I$, and let $B\in \Mat_{s,n}(K)$ 
be the coefficient matrix of $(\Lin_\M(f_1)_, \dots, \Lin_\M(f_s))$. 
Let $\sigma$ be a term ordering such that $\LT_\sigma(f_i)=z_i$ for $i=1,\dots,s$.
W.l.o.g., we may assume that $z_1 >_\sigma \cdots >_\sigma z_s$.
Then the submatrix of~$B$ consisting of the columns corresponding to~$Z$ 
is an upper triangular matrix $U_s \in \GL_s(K)$. By \cite{KLR1}, Proposition~1.6,  
we know that $\Lin_\M(f_i) \in \langle\Lin_\M(g_1), \dots, \Lin_\M(g_r)\rangle_K$ for $i=1,\dots, s$. 
Hence we deduce that there exists a matrix $C\in\Mat_{s,r}(K)$ such that $B = C \cdot A$, and therefore
$U_s = C\cdot A'$. Now that fact that $U_s$ is invertible implies the claims.
\end{proof}

Tuples~$Z$ for which the rank condition in this proposition holds are said to be 
{\bf of top rank} with respect to $\Lin_\M(I)$. They will turn out to be quite useful later.
Notice that the property of~$Z$ to be of top rank does not depend on the choice of a system of 
generators of $\Lin_\M(I)$.

\medskip
Particularly desirable re-embeddings are defined as follows.

\begin{definition}\label{def-optimal}
Let $I\subseteq \M$ be an ideal in~$P$. 
\begin{enumerate}
\item[(a)] If there exist a polynomial ring $P'=K[y_1,\dots,y_m]$,  
an ideal~$I'$ in~$P'$, and  an isomorphism
$\Phi:P/I \to P'/I'$ of $K$-algebras such that~$m$ is the minimum possible number, 
we say that $\Phi$ is an {\bf optimal re-embedding} of $I$ (or of $P/I$).

\item[(b)] An optimal re-embedding of~$I$ which is obtained from a $Z$-separating re-embedding 
will be called an {\bf optimal separating re-embedding} of~$I$.

\end{enumerate}
\end{definition}

Not every optimal re-embedding of~$I$ is an optimal separating re-embedding, as for instance
Examples~\ref{ex-optimalZ} and~\ref{ex-crachiola}.b show.
On the positive side, in~\cite{KLR1}, Corollary~4.2, we proved the following criterion 
for a $Z$-separating re-em\-bed\-ding to be optimal. 

\begin{proposition}\label{prop-Z-optimal}
Let $I\subseteq \M$ be an ideal in~$P$, let $Z$ be a tuple of indeterminates,
and let~$F$ be a coherently $Z$-separating tuple of polynomials in~$I$. 
If~$\# Z = \dim_K(\Lin_\M(I))$, then the 
corresponding $Z$-separating re-embedding~$\Phi$ of~$I$ is optimal.
\end{proposition}

For instance, in Example~\ref{ex-ElimBySubst}, we have $\# Z = 2 = \dim_\QQ(\Lin_\M(I))$, and therefore
the isomorphism $\Phi: P/I \longrightarrow \QQ[a,w] / \langle a^5 w^2 - a^2 w^2 \rangle$ 
is an optimal separating re-embedding. However, condition  $\# Z = \dim_K(\Lin_\M(I))$ is not necessary 
as the following example shows (see~\cite{KLR2}, Example~6.6).

\begin{example}\label{ex-optimalZ}
Let $P =\QQ[x, y]$, and let $I = \langle x - y^2,\, y - x\rangle$. In this case
we have $\M =\langle x, y \rangle$ and $\Lin_\M(I) = \langle x, y \rangle$.
For the tuple $Z=(x)$, we have $s=1$ and the inequality 
$\dim_{\QQ}(\Lin_\M(I)) =2 > s=1$. Nevertheless, the isomorphism
$\Phi: P/I \To \mathbb{Q}[y]/\langle y-y^2 \rangle$ is an optimal re-embedding,
since $P/I$ cannot be isomorphic to~$\QQ$, as~$\mathbb{Q}[y]/\langle y-y^2\rangle$ 
is not an integral domain.
\end{example}

In~\cite{KLR1}, Definition~3.3.d, we used the  name  ``optimal separating re-embedding" 
for $Z$-separating re-embeddings such that $\#Z$ is as large as possible. These will be called 
{\it best}\/ separating re-embeddings here (see Definition~\ref{def-best}), as we now feel 
that an optimal separating re-embedding should be {\it truly}\/ optimal.

\bigskip\bigbreak
%
%

\section{Re-Embeddings of Non-Negatively Graded Ideals}
\label{Re-Embeddings of Non-Negatively Graded Ideals}

It is time to move to the first main topic of this paper which is to find
homogeneous separating tuples of polynomials and the corresponding
re-embeddings in the non-negatively graded case.

Let us begin by fixing the setting. We let~$K$ be a field, $P=K[x_1,\dots,x_n]$ a polynomial
ring, $X=(x_1,\dots,x_n)$ its tuple of indeterminates, and $\M=\langle x_1,\dots,x_n\rangle$.
We use several kinds of gradings on~$P$. 
For the standard grading given by $\deg(x_1) = \cdots = \deg(x_n) = 1$, 
we denote the degree of a polynomial~$f$ by~$\deg(f)$.
For an arbitrary $\ZZ$-grading given by a matrix 
$W=(w_1,\dots,w_n)\in\Mat_{1,n}(\ZZ)$, we use $\deg_W(g)$ to 
denote the $W\!$-degree of a polynomial~$g$. 
Recall that this grading is given by $\deg_W(x_i)=w_i$ for $i=1,\dots,n$ (see~\cite{KR2}, Chapter~4).
Moreover, the grading by~$W$ is called {\bf non-negative} if $w_i\ge 0$ for $i=1,\dots,n$,
and it is called {\bf positive} if $w_i >0$ for $i=1,\dots,n$. In both cases the non-zero
elements of~$K$ have degree~0.  The standard grading is a special case of
a positive grading.

In this setting we assume that we are given a $W\!$-homogeneous
ideal~$I$ of~$P$, and that~$I$ is generated by elements of positive degree,
i.e., that $I\cap P_0 = \{0\}$. Notice that, if the grading given by~$W$ is positive, 
this condition simply means that the ideal is proper. If the grading is non-negative, 
the condition implies that~$I$ is contained in the ideal generated by the indeterminates 
of positive $W$-degree, hence that $I\subseteq \M$, and that the grading is non-trivial.

Moreover, we write $I=\bigoplus_{d\ge 1} I_d$ for the
decomposition of~$I$ into $W\!$-homogeneous components, we
let $I_{\langle d\rangle}$ be the ideal generated by~$I_d$, and we let~$R$ be the $W\!$-graded 
$K$-algebra $R=P/I$.

\begin{remark}{\bf ($W\!$-Homogeneous $Z$-Separating Tuples)}\label{rem-linzero}\\
Assume that we are in the above setting. Let $I$ be an ideal in $P$ generated by 
elements of positive degree, and let $z$ be an indeterminate in~$X$.
\begin{enumerate}
\item[(a)] If $f\in I$ is a $z$-separating polynomial, then the homogeneous 
component~$f_d$ of~$f$ of degree $d=\deg_W(z)$ is also contained in~$I$ 
and is a $z$-separating polynomial, as well.
Therefore, when needed, we may assume that all $Z$-separating tuples in~$I$
consist of $W\!$-homogeneous polynomials.

\item[(b)] If the indeterminate~$z$ is of $W\!$-degree zero, it is 
not contained in $\Supp(\ell)$ for any $\ell\in\Lin_\M(I)$. 
In particular, the ideal~$I$ contains no $z$-separating $W\!$-homogeneous polynomial.
The reason is that if $z\in \Supp(\Lin_\M(f))$ for some $W\!$-ho\-mo\-ge\-neous 
polynomial $f\in I$, then~$f$ is $W\!$-ho\-mo\-ge\-neous of
degree zero, and thus $f\in P_0$ in contradiction to $I\cap P_0= \{0\}$.

\end{enumerate}
\end{remark}

The next theorem allows us to check effectively whether a $Z$-separating tuple
exists in~$I$ for a given tuple of indeterminates~$Z$. To simplify its presentation,
we introduce the following abbreviations.

\begin{definition}\label{def-linZ}
Let $Z=(z_1,\dots,z_s)$ be a tuple of distinct indeterminates in $X=(x_1,\dots,x_n)$
which have positive $W\!$-degrees.
\begin{enumerate}
\item[(a)] Given a homogeneous polynomial $g\in P$, we write it in the form
$g = c_1 z_1 + \cdots + c_s z_s + h$ where $c_1,\dots,c_s \in P_0$ and where
the support of $h\in P$ does not contain any term of the form $t\, z_i$ with $t\in P_0$.

Then we call $\lin_Z(g) = c_1 z_1 +\cdots + c_s z_s$ the
{\bf $Z$-linear part} of~$g$ and $\cvec_Z(g) = (c_1,\dots,c_s) \in P_0^s$
the {\bf $Z$-coefficient vector} of~$g$. 

\item[(b)] Now let $G=(g_1,\dots,g_r)$ be a tuple of $W\!$-homogeneous polynomials.
Then the tuple $\lin_Z(G) \;=\; (\lin_Z(g_1),\dots, \lin_Z(g_r))$
is called {\bf $Z$-linear part} of~$G$.

Moreover, the matrix $\cmat_Z(G) \in \Mat_{s,r}(P_0)$ whose columns are the
$Z$-coefficient vectors $\cvec(g_1)^{\rm tr}, \dots, \cvec(g_r)^{\rm tr}$
is called the {\bf $Z$-coefficient matrix} of~$G$.
In particular, we have $Z\cdot \cmat_Z(G) =  \lin_Z(G)$.
\end{enumerate}
\end{definition}

Now we are ready to formulate and prove the main theorem of this section.
In the following we let $\{e_1,\dots,e_s\}$ denote the standard basis of $P_0^s$, i.e.,
$e_i = (0,\dots,0,1,0,\dots,0) \in P_0^s$, where~$1$ occurs in the $i$-th position.

\begin{theorem}{\bf  (Existence of Homogeneous $Z$-Separating Tuples)}\label{thm-unitvec}\\
Let $P=K[x_1,\dots,x_n]$ be a polynomial ring which is non-negatively $\ZZ$-graded by
$W \in\Mat_{1,n}(\ZZ)$, let~$I$ be a $W\!$-homogeneous ideal in~$P$ such that
$I \cap P_0 = \{0\}$, let $G=(g_1,\dots,g_r)$ be a tuple of $W\!$-homogeneous polynomials 
which generates~$I$, and let $Z =(z_1,\dots,z_s)$ be a tuple of distinct indeterminates in~$X$
of positive degrees. Then the following conditions are equivalent.
\begin{enumerate}
\item[(a)] There exists a $Z$-separating tuple of polynomials in~$I$.

\item[(b)] There exists a coherently $Z$-separating tuple of polynomials in~$I$.

\item[(c)] The elements $z_1,\dots, z_s$ are contained in the $P_0$-module generated 
by the set $\{\lin_Z(g_1), \dots, \lin_Z(g_r)\}$.

\item[(d)] The vectors  $e_1, \dots, e_s$ are contained in the
$P_0$-module generated by the set $\{\cvec_Z(g_1), \dots, \cvec_Z(g_r)\}$.

\end{enumerate} 
\end{theorem}

\begin{proof}
It follows from Remark~\ref{rem-sepandcoherentlysep} that claims~(a) and~(b) are equivalent.

Next we show that~(b) implies~(c).
By assumption and Remark~\ref{rem-linzero}.a, there exists a tuple $F = (f_1,\dots,f_s)$ 
of $W\!$-homogeneous polynomials in~$I$ such that $z_i= \lin_Z(f_i)$ for $i=1,\dots,s$.
We write $f_i = p_{i1} g_1 + \cdots + p_{ir} g_r$ with homogeneous polynomials $p_{ij}\in P$
and $g_j = c_{j1} z_1 + \cdots + c_{js} z_s + h_j$, where $c_{jk}\in P_0$ and $h_j\in P$ 
contains no term of the form $tz_k$ with $t\in P_0$ and $k\in \{1,\dots,s\}$. Then we examine the equalities
$$
z_i \;=\; \lin_Z(f_i) \;=\; \lin_Z \big( p_{i1} (c_{11} z_1 + \cdots + c_{is} z_s + h_1) + \cdots +
p_{ir} (c_{r1} z_1 + \cdots + c_{rs} z_s + h_r) \big)
$$
and see that $z_i \in \langle  \lin_Z(g_1), \dots, \lin_Z(g_r)\rangle_{P_0}$, as claimed.

To prove that~(c) implies~(d), it suffices to note that $\{z_1,\dots,z_s\}$ is a basis
of the $P_0$-module $P_0 z_1 + \cdots + P_0 z_s$ and two write down the coordinate tuples.

Lastly, we prove that~(d) implies~(a). By assumption, there exits a matrix
$A\in \Mat_{r,s}(P_0)$ such that $I_s = \cmat_Z(G)\cdot A$.
This implies $Z = Z  \cdot \cmat_Z(G)\cdot A = \lin_Z(G) \cdot A$. Consequently,
the polynomials in $F =(f_1,\dots,f_s) = \lin_Z(G) \cdot A$ satisfy $\lin_Z(f_i)=z_i$
for $i=1,\dots,s$. Now we construct a term ordering $\sigma$ on $\mathbb T^n$ as follows. 
We distinguish two cases:

{\it Case 1:} If there is only one indeterminate of positive degree in~$Z$, we put
$(1,0,\dots,0)$ into the first row of a matrix and extend it to a term ordering matrix.

{\it Case 2:} If there are at least two indeterminates in~$Z$, 
let $(d_1, \dots, d_s)$ be the tuple of their $W\!$-degrees. W.l.o.g., assume that this tuple  
is ordered by non-decreasing $W\!$-degree. Then construct a matrix whose first row
is  $(d_1,\dots,d_s, 0,\dots,0)$, the second row is $(d_1-1,\dots, d_s-1,0,\dots,0)$, 
and extend everything to a term ordering matrix.

Using this term ordering, it remains to verify that $\LT_\sigma(f_i)=z_i$
for $i=1,\dots,s$. By the construction of the matrix, we know that~$z_i$ is larger w.r.t.~$\sigma$ 
than any term of $W\!$-degree~$d_i$ which is a product of indeterminates, 
except for terms of the form $tz_i$ with $t\in P_0\setminus K$. 
Since the latter terms do not occur in $\lin_Z(f_i)$, and since~$f_i$ is $W\!$-homogeneous 
of degree~$d_i$, we get $\LT_\sigma(f_i)=z_i$, as desired.
Altogether, the proof of the theorem is complete.
\end{proof}

Based on this theorem and its proof, we obtain the following
algorithm which, given a tuple~$Z$ of distinct indeterminates,
checks if $Z$ is separating for~$I$ and, if it is, computes a $Z$-separating 
tuple of polynomials in~$I$. Notice that it is essential to assume 
that the grading given by~$W$ is non-negative.

\begin{algorithm}{\bf (Finding Homogeneous  $Z$-Separating Tuples)}\label{alg-findF}\\
Let $P=K[x_1,\dots,x_n]$ be a polynomial ring which is non-negatively $\ZZ$-graded by
$W \in\Mat_{1,n}(\ZZ)$, let~$I$ be a $W\!$-homogeneous ideal in~$P$ such that
$I \cap P_0 = \{0\}$, let $G=(g_1,\dots,g_r)$ be a tuple of $W\!$-homogeneous polynomials 
which generates~$I$, and let $Z =(z_1,\dots,z_s)$ be a tuple of distinct indeterminates in~$P$
of positive degrees.  
Consider the following instructions.

\begin{enumerate}
\item[(1)] For $i=1,\dots,r$, find $\cvec_Z(g_i)$.

\item[(2)] 
If the $P_0$-module membership problem to find $(a_{i1},\dots,a_{ir}) \in P_0^r$
such that $e_i = a_{i1} \cvec_Z(g_i) + \cdots + a_{ir} \cvec_Z(g_r)$ 
has no solution for some $i\in\{ 1,\dots,s\}$, return {\tt "Z is not separating for I"} and stop.

\item[(3)] For $i=1,\dots, s$, let $(a_{i1},\dots,a_{ir}) \in P^r$ be a solution of the $P_0$-module
membership problem in~(2) and form the polynomial $f_i = a_{i1} g_1 + \cdots + a_{ir} g_r$.

\item[(4)] Return the tuple $F=(f_1, \dots, f_s)$ and stop.
\end{enumerate}

This is an algorithm which checks if  a homogeneous $Z$-separating tuple
of polynomials exists in~$I$, and in that case computes such a tuple~$F$.

Moreover, for each degree~$d$, let $Z_d$ be the tuple of indeterminates of degree~$d$
in~$Z$, and let~$F_d$ be the corresponding subtuple of~$F$. Then~$F_d$ is a coherently
$Z_d$-sep\-a\-ra\-ting tuple.
\end{algorithm}Theorem

\begin{proof}
The first claim follows immediately from Theorem~\ref{thm-unitvec} and its proof.
To verify the additional claim, note that step~(2) implies $\lin_{Z_d}(f_i)= z_i$ for every indeterminate~$z_i$
in~$Z_d$.
Since we have $\deg_W(f_i)=d$, the indeterminates in~$Z_d$ divide only the $Z$-linear part of~$f_i$,
and therefore the only term in $\Supp(f_i)$ divisible by some indeterminate~$z_j$ in~$Z_d$ is~$z_i$.
Hence the tuple $F_d$ is coherently $Z_d$-sep\-a\-ra\-ting. 
\end{proof}

The following example illustrates Theorem~\ref{thm-unitvec} and the above algorithm.

\begin{example}\label{ex-cvec}
Let $P = \QQ[a,b,x,y,z]$ be graded by $W = (0,0,1,1,1)$, and let 
$I = \langle g_1, g_2, g_3 \rangle$, where 
$g_1 = x - a^3y +z$,   $g_2 = x - (ab+1)y$, and $g_3 = a^2y +az$.
The ideal~$I$ is $W\!$-homogeneous and satisfies $\Lin_\M(I) = \langle x, y, z\rangle_\QQ$
as well as $I\cap P_0 = I\cap\QQ[a,b] = \{ 0 \}$.

For $Z =(x, y)$, we have $\cvec_Z(g_1) = (1, -a^3)$,
$\cvec_Z(g_2) = (1, -ab -1)$, and $\cvec_Z(g_3) = (0, a^2)$. 
When we try to represent~$e_1$ and~$e_2$ using these coefficient vectors, we get
\begin{align*}
e_1 &\;=\; 1\cdot \cvec_Z(g_1) + a\cdot\cvec_Z(g_3)\\
e_2 &\;=\; (-ab +1)\cdot \cvec_Z(g_1) + ( ab -1)\cdot \cvec_Z(g_2) + (-a^2b +b^2 +a)\cdot \cvec(g_3)
\end{align*}
This yields the homogeneous polynomials $f_1 = 1\cdot g_1 + a\cdot g_3 = x +z+a^2z$ and
$f_2 = (-ab {+}1)\cdot g_1  + (ab{-}1)\cdot g_2+ (-a^2b{+}b^2{+}a)\cdot g_3 = 
y + z -a^3bz +ab^2z +a^2z -abz$.
Altogether, the pair $(f_1,f_2)$ is a coherently $Z$-separating tuple.
\end{example}

As a natural continuation of the above algorithm, we now show how 
to pass from a homogeneous $Z$-separating tuple to a homogeneous coherently
$Z$-separating tuple.

\begin{algorithm}{\bf (Finding Homogeneous Coherently $Z$-Separating Tuples)\label{alg-findCohF}}\\
In the setting of Algorithm~\ref{alg-findF} assume that the output is a tuple~$F$.
W.l.o.g.\ let $\deg_W(z_1) \le \cdots \le \deg_W(z_s)$.
Consider the following instructions.
\begin{enumerate}
\item[(1)] Let $\tilde{f}_1 = f_1$.

\item[(2)] For $i=2,\dots,s$, rewrite the polynomial $f_i$ using the
tuple $(\tilde{f}_1,\dots,\tilde{f}_{i-1})$ and get a new polynomial $\tilde{f}_i$.

\item[(3)] Return the tuple $\widetilde{F} = (\tilde{f}_1,\dots, \tilde{f}_s)$ and stop.
\end{enumerate}
This is an algorithm which computes a homogeneous coherently $Z$-separating 
tuple~$\widetilde{F}$ of  polynomials  in~$I$.
\end{algorithm}

\begin{proof}
By Algorithm~\ref{alg-findF}, inside each degree~$d$ the tuple~$F$ is already coherently $Z_d$-separating.
Thus the only indeterminates~$z_j$ dividing a term in $\Supp(f_i)$ can be of lower
$W\!$-degree, and they are reduced in step~(2).
\end{proof}

Notice that in step~(2) it suffices to rewrite~$f_i$ using the polynomials $\tilde{f}_j$ 
of lower $W\!$-degree.
To finish this section we elaborate a special case of Algorithm~\ref{alg-findF}.
It allows us to weed out all indeterminates~$x_i$ for which no
$x_i$-separating homogeneous polynomial exists in~$I$ and will be a useful preprocessing step for  
further algorithms.

\begin{algorithm}{\bf (Finding All Separating Indeterminates)}\label{alg-FindSepPoly}\\
Let $P=K[x_1,\dots,x_n]$ be a polynomial ring which is non-negatively $\ZZ$-graded by
$W \in\Mat_{1,n}(\ZZ)$, let $g_1,\dots,g_r\in P$ be homogeneous of positive $W\!$-degree, and
let $I = \langle g_1,\dots,g_r\rangle$. Consider the following sequence of instructions.
\begin{enumerate}
\item[(1)] Let $Z$ be the empty tuple. For $i=1,\dots,n$, perform the Steps (2) and~(3).

\item[(2)] Let $\{g_{j_1},\dots,g_{j_t}\}$ be the set of all generators of~$I$
such that $\deg_W(g_{j_k}) = \deg_W(x_i)$. For $k=1,\dots,t$, we write 
$g_{j_k} = c_k \, x_i + h_k$ with $c_k\in P_0$ and $h_k\in P$ such that~$x_i$ divides no
term in the support of~$h_k$.

\item[(3)] If the ideal membership problem  $1\in \langle c_1,\dots,c_t\rangle_{P_0}$ in~$P_0$ 
has a solution, append~$x_i$ to~$Z$.

\item[(4)] Return the resulting tuple~$Z$ and stop.
\end{enumerate}
This is an algorithm which computes the tuple~$Z$ of all indeterminates~$x_i$ in~$X$ such that
a $x_i$-separating homogeneous polynomial exists in~$I$, i.e., the tuple of all separating
indeterminates for~$I$.
\end{algorithm}

\bigskip\bigbreak
%
%

\section{Best Separating Re-Embeddings}
\label{sec4}

As before, we let $P=K[x_1,\dots,x_n]$ be non-negatively graded by~$W\!$, and
we let~$I$ be a $W\!$-homogeneous ideal in~$P$ such that $I\cap P_0 = \{0\}$.
Our goal in this section is to find a {\it best} separating re-embedding of~$I$ in the following sense.

\begin{definition}{\bf (Best Separating Re-Embeddings)}\label{def-best}\\
Let $I\subseteq \M$ be an ideal, and let $s \le \dim_K(\Lin_\M(I))$ be the largest number 
such that there exist a tuple of indeterminates $Z=(z_1,\dots,z_s)$ and a coherently $Z$-separating tuple~$F$
of polynomials in~$I$. 

Then such a tuple~$Z$ is called a {\bf best tuple of separating indeterminates} for~$I$ and
the corresponding $Z$-separating re-embedding $\Phi:\; P/I \longrightarrow K[X\setminus Z] / (I\cap K[X\setminus Z])$
is called a {\bf best separating re-embedding} of~$I$. 
\end{definition}

By Proposition~\ref{prop-Z-optimal}, we know that if $\#Z = \dim_K(\Lin_\M(I))$ then the corresponding
best separating re-embedding of~$I$ is actually optimal. However, in the last part of Section~\ref{sec2}
we saw that if $\#Z < \dim_K(\Lin_\M(I))$ then the corresponding re-embedding can be optimal or not. 

In order to discover a best separating re-embedding in the case where all generators of~$I$
have the same degree~$d$, we can form the set of all candidate indeterminates
and search through its subsets as follows.

\begin{algorithm}{\bf (A Best Tuple of Separating Indeterminates in Degree~$d$)}\label{alg-BestZSepInDegd}\\
Let $d\ge 1$, and let~$I$ be a $W\!$-homogeneous ideal in~$P$ such that $I \cap P_0 = \{0\}$ and such that~$I$
is generated by homogeneous polynomials $G = (g_1, \dots, g_r)$ of $W\!$-degree~$d$.
Consider the following sequence of instructions.
\begin{enumerate}
\item[(1)] Using Algorithm~\ref{alg-FindSepPoly}, calculate the set $S = \{ z_1,\dots,z_u\}$
consisting of all indeterminates~$z_i$ of $W\!$-degree~$d$ for which~$I$ contains
a $z_i$-separating polynomial.
If $S {=} \emptyset$, return {\tt "No separating indeterminates in degree d"} and stop.

\item[(2)] Let $m=\min\{s,\, \dim_K(\Lin_\M(I)) \}$. Starting with $m'=m$, 
form all $m'$-element subsets~$S'$ of~$S$ and execute step~(3) for each of them 
until it returns a result. If it returns no result for any subset~$S'$ 
of a given cardinality~$m'$, reduce~$m'$ by one and continue.

\item[(3)] Let $Z=(z_1,\dots,z_s)$ be the tuple of all indeterminates in~$S'$, ordered
arbitrarily. Check if~$Z$ is of top rank. If this is the case, use 
Algorithm~\ref{alg-findF} to check if there exists a $Z$-separating 
tuple of $W\!$-homogeneous polynomials in~$I$. In this case
compute such a tuple~$F$. Otherwise, continue with the next set~$S'$ in step~(2). 

\item[(4)] Return the pair $(Z,F)$ and stop.
\end{enumerate} 
This is an algorithm which computes a pair $(Z,F)$, where~$Z$ is a best separating tuple 
of indeterminates of $W\!$-degree~$d$ and where~$F$
is a $Z$-separating tuple of homogeneous polynomials in~$I$.
\end{algorithm}

\begin{proof}
Finiteness of the algorithm follows if we show that the loop in steps~(2) and~(3)
is finite. This is a consequence of the hypothesis that~$S$ is finite, as it consists of
indeterminates, and thus step~(1) implies that the loop stops at the very latest 
when $m'=1$.

To prove correctness, we note that the loop in steps~(2) and~(3) works through 
all subsets~$S'$  of~$S$ in descending order of cardinality. Therefore the tuple~$Z$
found eventually has the maximal possible cardinality and is hence a best tuple
of separating indeterminates for~$I$.
\end{proof}

A simple modification of this algorithm finds all best coherently separating
tuples in one degree as follows.

\begin{remark}{\bf (All Best Separating Tuples in Degree~$d$)}\label{rem-allbest}\\
Algorithm~\ref{alg-BestZSepInDegd} can be slightly modified to compute not only the first 
tuple~$Z$ which yields a best separating re-embedding of~$I$, but all such tuples.  

Let $d\ge 1$, and execute Algorithm~\ref{alg-BestZSepInDegd} until it finds one number~$m'$
and one set $S'\subseteq S$ of cardinality~$m'$ which yields a separating tuple of indeterminates~$Z$ 
in step~(3). Then form all other $m'$-element subsets~$S'$
of~$S$ and perform step~(3) for each of them. Collect all resulting pairs $(Z,F)$ in a set~$B$.

Then~$B$ contains all tuples~$Z$ of indeterminates of $W\!$-degree~$d$, each together with
a $Z$-separating tuple of polynomials~$F_Z$, such that the re-embeddings defined by the various tuples~$F_Z$
are {\it all}\/ best separating re-embeddings of~$I$.
\end{remark}

Here is a concrete example for the algorithm described in this remark.

\begin{example}\label{ex-allbest}
Let $P = \QQ[a,b,v,x,y,z,w]$ be graded by $W = (0,0,1,2,2,2,2)$, and let 
$I = \langle g_1, g_2, g_3, g_4 \rangle$, where 
$g_1 = x -z -w + ax + b v^2 - v^2$,  $g_2 = 3bx + y - v^2$, $g_3 = abx$, and $g_4 = ax - v^2$.
The ideal~$I$ is $W\!$-homogeneous and satisfies $\Lin_\M(I) = \langle x, y, z, w\rangle_\QQ$
as well as $I\cap P_0 = I\cap\QQ[a,b] = \{ 0 \}$.

Since all generators of~$I$ have $W$-degree 2, we consider $d=2$. Executing
Algorithm~\ref{alg-BestZSepInDegd} yields that the choice of $S'=\{y,x\}$
in step~(3) succeeds. Thus the result is the pair
$((y,x), (y -3b^2 v^2 +3bz +3bw - v^2,\; x + b v^2 -z -w))$.

The modified algorithm continues to check 2-element subsets~$S'$
of~$S$ and returns the following pairs.
\begin{align*}
((y,x), \;& (y - 3b^2 v^2 + 3bz + 3bw - v^2,\; x + b v^2 -z -w))\cr
((y,z), \;&( y + 3bx - v^2,\;  z - b v^2 -ax + v^2 -x +w))\cr
((y,w), \;& (y + 3bx - v^2,\;  w - b v^2 -ax + v^2 -x +z))
\end{align*}
Thus there exist three different best separating re-embeddings of~$I$. 
\end{example}

Our main goal is to compute a best separating re-embedding of~$I$ 
for an arbitrary $W\!$-homogeneous ideal~$I$ such that $I\cap P_0 = \{0\}$. 
For this, our strategy is to first show that we can find~$Z$ degree-by-degree
and then use Algorithm~\ref{alg-BestZSepInDegd} in each degree.

The next proposition says that this strategy is correct.
To ease the notation, for a given set or tuple~$S$ of $W\!$-homogeneous polynomials
and a degree $d\ge 0$, we let $S_d$ be the subset or subtuple consisting of all elements
of~$S$ of $W\!$-degree~$d$.

\begin{proposition}{\bf (Homogeneous Separating Tuples Degree by Degree)}\label{prop-deg-by-deg}\\
Let $G=(g_1,\dots,g_r)$ be a tuple of homogeneous polynomials of positive $W\!$-degrees
which generates an ideal~$I$ in~$P$. For every $d\ge 1$, let $I_d = \langle G_d\rangle$.
\begin{enumerate}
\item[(a)] A tuple of homogeneous polynomials~$F$ in~$I$ is $Z$-separating if and only if 
for every $d\ge 1$ the tuple $F_d$ is $Z_d$-separating for the ideal~$I_d$. 

\item[(b)] A tuple of homogeneous polynomials~$F$ in~$I$ is coherently $Z$-separating if and only if 
for every $d\ge 1$ the tuple $F_d$ is coherently $Z_d$-separating for the ideal~$I_d$
and $F_d$ is reduced w.r.t.\ all tuples $F_{d'}$ with $d'<d$.

\item[(c)] A tuple~$Z$ is a best tuple of separating indeterminates for~$I$ if and only 
if for every $d\ge 1$ the tuple~$Z_d$ is a best tuple of separating indeterminates for~$I_d$.

\end{enumerate}
\end{proposition}

\begin{proof}
Claim~(a) follows from Theorem~\ref{thm-unitvec}. Claim~(b) is a consequence of~(a)
and Algorithm~\ref{alg-findCohF}. Finally, claim~(c) follows from~(a) and the definition
of a best tuple of separating indeterminates.
\end{proof}

By combining Algorithm~\ref{alg-BestZSepInDegd} with this proposition, 
we arrive at the following general method for computing a best tuple of separating indeterminates 
in the non-negatively graded case.

\begin{algorithm}{\bf (Computing a Best Tuple of Separating Indeterminates)}\label{alg-BestPossible}\\
Let $G=(g_1,\dots,g_r)$ be a tuple of $W\!$-homogeneous polynomials such that
the ideal $I=\langle g_1,\dots,g_r\rangle$ satisfies $I \cap P_0 = \{0\}$.
Consider the following steps.
\begin{enumerate}
\item[(1)] Using Algorithm~\ref{alg-FindSepPoly}, compute the set~$S$ of all
indeterminates~$x_i$ in~$X$ such that~$I$ contains an $x_i$-separating polynomial.

\item[(2)] Let $D = \{ \deg_W(x_i) \mid x_i \in S\}$. For every $d\in D$, 
let $I_d$ be the ideal generated by the $W\!$-homogeneous elements of degree~$d$ in~$G$.
Use Algorithm~\ref{alg-BestZSepInDegd} to compute a pair $(Z_d,F_d)$ such that~$Z_d$
is a best tuple of separating indeterminates of $W\!$-degree~$d$ and~$F_d$ is a $Z_d$-separating
tuple of polynomials.

\item[(3)] Return the pair $(Z,F)$ where $Z=\bigcup_{d\ge 1}$ and $F=\bigcup_{d\ge 1}$.
\end{enumerate}
This is an algorithm which computes a pair $(Z,F)$ such that~$Z$ is a best tuple of separating
indeterminates for~$I$ and~$F$ is a $Z$-separating tuple of homogeneous polynomials.
\end{algorithm}

If we are interested in a coherently $Z$-separating tuple, and thus a best separating
re-embedding of~$I$, it suffices to combine this algorithm with Proposition~\ref{prop-deg-by-deg}.d.
The following example shows the algorithm at work.

\begin{example}\label{ex-optimalgraded}
Let $P = \QQ[a, x, y, z, v, w]$ be graded by $W = (0,1,1,1,2,2)$, and let 
$I = \langle g_1, \dots g_6\rangle$, where $g_1 = x -ay$, $g_2 = y -az$, 
$g_3 = z -ax$,  $g_4 =  v -a^3w -axz$, $g_5 = v - w -aw +ayz$, and $g_6= a^4w$.
The ideal~$I$ is $W\!$-homogeneous, satisfies $\Lin_\M(I) = \langle x, y, z, v, w\rangle_\QQ$,
and we have $I\cap P_0 = I\cap\QQ[a] =\{0\}$. 

Our goal is to find a best separating re-embedding of~$I$.
We follow the steps of the algorithm.
\begin{enumerate}
\item[(1)] Using Algorithm~\ref{alg-FindSepPoly}, or by looking at $\Lin_\M(I)$ 
and the generators of~$I$, we see that it suffices to search for tuples~$Z$ 
whose entries are contained in $S = \{x, y, z, v, w\}$.

\item[(2)] We have $D=\{1,2\}$. In $W\!$-degree $d=1$, Algorithm~\ref{alg-BestZSepInDegd} 
returns the pair  $\big((x, y), (f_1, f_2)\big)$   where $f_1 = x-a^2z$ and $f_2 = y - az$.

\item[(2)] In $W\!$-degree $d=2$, Algorithm~\ref{alg-BestZSepInDegd} returns
$((v, w),\, (f_3,f_4))$, where  
$f_3 = v - a^4xz -a^4yz -axz$ and $f_4 = w -a^3xz -a^3yz +a^2xz +a^2yz -axz -ayz$.

\item[(3)] We get $Z= (x, y, v, w)$ and the $Z$-separating tuple $F=(f_1,f_2,f_3,f_4)$. 
\end{enumerate}

If we are looking for a coherently $Z$-separating tuple, we follow Algorithm~\ref{alg-findCohF}.
Thus we substitute $x \mapsto a^2z$ and $y \mapsto az$ in~$f_3$ and~$f_4$ and get 
$$
\widetilde{F} \;=\; (x -a^2z,\; \; y - az,\; \;  v-a^6z^2 -a^5z^2 -a^3z^2,\;\;  w -a^5z^2 -a^2z^2)
$$

The corresponding best separating re-embedding is the $\QQ$-algebra isomorphism
$\Phi:\; P/I \longrightarrow \QQ[a, z]/ (I\cap \QQ[a,y])$ which is induced by the $\QQ$-algebra 
homomorphism $\phi:\; P \longrightarrow \QQ[a, z]$ 
defined by $\phi(a) = a $, $\phi(x) = a^2z$, $\phi(y) = az$, $\phi(z) = z$,
$\phi(v) = a^6z^2 +a^5z^2 +a^3z^2$, and $\phi(w) = a^5z^2 +a^2z^2$.

As explained in Proposition~\ref{prop-ElimViaSubst}, the elimination ideal
$I\cap \QQ[a, z]$ is the result of performing these substitutions in the generators of~$I$.
Hence we get
$$
I\cap\QQ[a, z] \;=\; \langle 0,\; 0,\; -a^3z +z, \; a^8z^2 +a^6z^2,\;  0, \; a^9z^2 +a^6z^2\rangle \;=\; 
\langle (a^3-1)z, z^2 \rangle 
$$
Since $\dim_\QQ(\Lin_\M(I)) = 5$ and we found only a tuple~$Z$ with four
indeterminates, we cannot deduce from Proposition~\ref{prop-Z-optimal} that~$\Phi$ is an optimal 
re-embedding of~$I$.
However, it turns out that~$\Phi$ is indeed optimal. 
The reason is that $\dim(P/I) = \dim (\QQ[a, y]/ (I\cap \QQ[a, y]) =1$. 
To lower the number of indeterminates further, we would have to re-embed 
$I\cap \QQ[a,y]$ into a univariate polynomial ring. The only 1-dimensional 
residue class ring of a univariate polynomial ring is the ring itself. On the other hand, 
$\QQ[a,y]/ (I\cap \QQ[a,y])$ is not an integral domain. 

Later we will continue the discussion of this example (see Example~\ref{ex-fibers}).
\end{example}

The case of positively graded $K$-algebras is a special case of the theory developed so far.
The main difference is that $P_0=K$, and hence the vectors $\lin_Z(f)$ are simply $K$-linear 
combinations of the indeterminates in~$Z$.
Therefore this case can be treated using linear algebra techniques.
In particular, we automatically obtain tuples of top rank whose cardinality equals the dimension of 
$\Lin_\M(I),$ and hence provide optimal separating re-embeddings.
The following corollary, a version of which is given as Proposition~3.4 in~\cite{RT},
spells this out explicitly.

\begin{corollary}{\bf (Optimal Re-Embeddings of Positively Graded Ideals)}\label{cor-optimalgraded}\\ 
Let $P=K[x_1,\dots,x_n]$ be positively graded by~$W\in \Mat_{1,n}(\ZZ)$, let $G=(g_1,\dots,g_r)$
be a tuple of $W\!$-homogeneous polynomials in~$\M$, and let $I=\langle g_1,\dots,g_r\rangle$.
Consider the following sequence of instructions.
\begin{enumerate}
\item[(1)] Let $s=\dim_K(\Lin_\M(I))$. Choose a subtuple 
$(f_1,\dots,f_s)$ of~$G$ such that $(\Lin_\M(f_1),\dots,\Lin_\M(f_s))$
is a $K$-basis of $\Lin_\M(I)$.

\item[(2)] Choose a tuple of pairwise distinct indeterminates $Z=(z_1,\dots,z_s)$ 
which is of top rank and such that $z_i\in\Lin_\M(f_i)$ for $i=1,\dots,s$.
W.l.o.g.\ let $\deg_W(z_1) \le\cdots\le \deg_W(z_s)$. 

\item[(3)] By $K$-linearly interreducing $(f_1,\dots,f_s)$, 
compute  $W\!$-homogeneous polynomials $f_1',\dots,f_s'$ such that
$f_i' = z_i - h_i$ with $h_i\in P$ and $\Lin_\M(h_i) \in K[X\setminus Z]$.

\item[(4)] For $i=1,\dots,s-1$ and $j=i+1, \dots,s$, substitute $z_i$ in~$f_j'$ by~$h_i$.

\item[(6)] Let $\tilde{f}_1,\dots,\tilde{f}_s$ be the polynomials resulting from these substitutions.
Return the tuple $\widetilde{F} = (\tilde{f}_1,\dots,\tilde{f}_s)$ and stop.
\end{enumerate}
This is an algorithm which computes a coherently $Z$-separating tuple~$\widetilde{F}$ 
which consists of $W\!$-ho\-mo\-geneous polynomials in the ideal~$I$ and
which defines an optimal separating re-embedding 
$\Phi:\; P/I \cong K[X\setminus Z] / (I \cap K[X\setminus Z])$ of~$I$.
\end{corollary}

\begin{proof}
Step~(1) can be executed because the linear parts of~$g_1,\dots,g_r$ generate $\Lin_\M(I)$.
Step~(2) makes sure that the $K$-vector space $\langle f_1,\dots,f_s\rangle_K$
contains polynomials with linear parts as required by step~(3).
Then step~(3) performs the linear algebra operations of Algorithm~\ref{alg-BestZSepInDegd}.
The homogeneity of the polynomials~$f_i'$ implies that~$z_i$ does not divide any element 
in $\Supp(h_i)$. Hence $f_i'$ is $z_i$-separating for $i=1,\dots, s$.

Next, steps~(4) and~(5) perform both the interreductions required by 
Algorithm~\ref{alg-BestZSepInDegd} and the rewriting required by Proposition~\ref{prop-deg-by-deg}.b.
The final tuple has $s=\dim_K(\Lin_\M(I))$ elements. Thus the corresponding $Z$-separating
re-embedding of~$I$ is an optimal separating re-embedding.
\end{proof}

Examples for the application of this algorithm are contained in the next sections.

\bigskip\bigbreak
%
%

\section{Optimal Re-Embeddings and Freeness of the Fibers}
\label{sec5}

As before, we consider a polynomial ring~$P$ in~$n$ indeterminates over a field~$K$ which
is non-negatively graded by $W=(w_1,\dots,w_n) \in \Mat_{1,n}(\ZZ)$. However, now we assume that
$w_1=\cdots = w_m=0$ and $w_{m+1},\dots,w_n >0$ for some $m\in \{0,\dots,n-1\}$.
Then we denote the tuple of the first~$m$ indeterminates by $\Xz=(a_1,\dots,a_m)$, 
so that $P_0 = K[a_1,\dots,a_m]$, and the tuple of the remaining indeterminates
by $\Xp=(x_{m+1},\dots,x_n)$. Furthermore, we let $L=K(a_1,\dots,a_m)$.

Given a $W\!$-homogeneous ideal~$I$ in~$P$ such that $I \cap P_0 = \{0\}$, we can view
the injective $K$-algebra homomorphism  $P_0 \longrightarrow P/I$ geometrically as a family
of schemes, parametrized by the affine space $\Spec(P_0)$. 

In this section we study optimal re-embeddings of the ideals defining these fibers and
describe when their coordinate rings are isomorphic to polynomial rings. In the next section we
will then study the re-embeddings and the freeness of the entire family.

\begin{definition}\label{def-fibers}
Let $I$ be a $W\!$-homogeneous ideal in~$P$ with $I\cap P_0 = \{0\}$, and let~$\phi$
be the canonical $K$-algebra homomorphism $\phi:\; P_0 \longrightarrow P/I$. 
\begin{enumerate}
\item[(a)] The ideal $I_L = I\, L[\Xp]$ is called the {\bf generic fiber ideal}, and
the induced $L$-algebra homomorphism $\phi_L:\; L \longrightarrow L[\Xp] / I_L$ (resp.\ 
the $L$-algebra $L[\Xp]/I_L$) is called the {\bf generic fiber} of~$\phi$.

\item[(b)] For $\Gamma = (c_1,\dots,c_m) \in K^m$, let $\m_\Gamma
= \langle a_1-c_1, \dots, a_m-c_m \rangle \subseteq P_0$.
Then the ideal $I_\Gamma = I \cdot (P_0 / \m_\Gamma)[\Xp]$ is called the
{\bf (special) fiber ideal} of~$\phi$ over~$\Gamma$.
The induced $K$-algebra homomorphism 
$\phi_\Gamma:\; K \longrightarrow K[\Xp] / I_\Gamma$ (resp.\ the $K$-algebra
$K[\Xp]/I_\Gamma$) is called the {\bf (special) fiber} of~$\phi$ over~$\Gamma$.

\end{enumerate}
\end{definition}

\begin{remark}\label{rem-surjective}
Notice that in part~(b) of this definition we have $P_0/\m_\Gamma \cong K$ and
we obtain $I_\Gamma \subseteq \langle X^+ \rangle \subseteq K[\Xp]$.
Consequently, the ideal~$I_\Gamma$ is a proper ideal in $K[\Xp]$ for every~$\Gamma$.
Hence the map $\Spec(P/I) \To \Spec(P_0)$ is surjective.
Since the grading is non-negative, the individual fibers of the corresponding morphism 
$\Proj(P/I) \longrightarrow \Spec(P_0)$ can be viewed as subschemes of the weighted 
projective space $\Proj(K[\Xp])$ (see~\cite{BR}).
\end{remark}

The following proposition implies that, under the assumption $I\cap P_0 = \{0\}$,
a coherently $Z$-separating tuple of $W\!$-homogeneous polynomials in~$I$
maps to a coherently $Z$-separating tuple in all fiber ideals.

\begin{proposition}{\bf (Coherently Separating Tuples in Fiber Ideals)}\label{prop-CohSepInFibers}\\
Let~$I$ be a $W\!$-homogeneous ideal in~$P$ with $I\cap P_0 = \{0\}$, and let
$Z=(z_1,\dots,z_s)$ be a tuple of indeterminates such that there is
a coherently $Z$-separating tuple of $W\!$-homogeneous polynomials $F=(z_1-h_1, \dots, z_s-h_s)$ 
in~$I$.
\begin{enumerate}
\item[(a)] We have $Z\subseteq \Xp$ and $h_i \in K[\Xz,\Xp{\setminus} Z]$ for $i=1,\dots, s$. 

\item[(b)] The tuple $F$ is a coherently $Z$-separating tuple in the generic fiber ideal~$I_L$.

\item[(c)] For $\Gamma = (c_1,\dots,c_m) \in K^m$, the tuple $F_\Gamma$ obtained from~$F$ by
the substitutions $a_i \mapsto c_i$ for $i=1,\dots,m$, is coherently $Z$-separating for the
fiber ideal~$I_\Gamma$.

\end{enumerate}
\end{proposition}

\begin{proof}
Claim (a) follows from Remark~\ref{rem-linzero}.b.
To show~(b) and~(c), it suffices to note that in both cases the indeterminates
in~$\Xz$ are transformed into elements of the base field, and therefore
the polynomials in~$F$ and~$F_\Gamma$, respectively, remain homogeneous and coherently $Z$-separating.
\end{proof}

Fiber ideals are homogeneous ideals contained in $L[\Xp]$ and $K[\Xp]$,
respectively, and therefore positively graded. This allows us to apply
Corollary~\ref{cor-optimalgraded} and get optimal re-embeddings of the fibers
as follows.

\begin{proposition}{\bf (Optimal Re-Embeddings of Fiber Ideals)}\label{prop-OptFibers}\\
Let~$I$ be a $W\!$-homogeneous ideal in~$P$ such that $I\cap P_0 =\{0\}$. 
\begin{enumerate}
\item[(a)] There exists a tuple of indeterminates~$Z$ 
in~$\Xp$ which yields an optimal separating re-em\-bed\-ding 
$$
\Phi_L:\; L[\Xp] / I_L \;\longrightarrow\; L[\Xp{\setminus} Z] / (I_L \cap L[\Xp{\setminus} Z])
$$
of the generic fiber ideal~$I_L$.

\item[(b)] Let $\Gamma=(c_1,\dots,c_m) \in K^m$. There exists a tuple of indeterminates~$Z$
in~$\Xp$ which yields an optimal separating re-em\-bed\-ding 
$$
\Phi_\Gamma:\; K[\Xp] / I_\Gamma  \;\longrightarrow\;  K[\Xp{\setminus} Z] / (I_\Gamma \cap K[\Xp{\setminus} Z])
$$
of the special fiber ideal~$I_\Gamma$.

\end{enumerate}
\end{proposition}

\begin{proof}
Note that~$I_L$ and~$I_\Gamma$ are proper homogeneous ideals in positively graded rings. 
Hence the claim follows by applying Corollary~\ref{cor-optimalgraded}. 
\end{proof}

The following example is a re-interpretation and extension of Example~\ref{ex-optimalgraded}
and  illustrates this proposition.

\begin{example}\label{ex-fibers}
As in Example~\ref{ex-optimalgraded}, let $P = \QQ[a, x, y, z, v, w]$ be graded by $W = (0,1,1,1,2,2)$, 
and let $I = \langle g_1, \dots, g_6\rangle$, where $g_1 = x -ay$, $g_2 = y -az$, 
$g_3 = z -ax$,  $g_4 =  v -a^3w -axz$, $g_5 = v - w -aw +ayz$, and $g_6= a^4w$.

Previously we computed an optimal re-embedding of the ideal~$I$ and found
$\Phi:\; P/I \longrightarrow \QQ[a, z]/(I\cap \QQ[a, z])$, where 
$$
I\cap \QQ[a, z] \;=\; \langle (a^3-1)z, z^2\rangle 
\;=\;  \langle a^3-1, z^2 \rangle  \cap \langle z\rangle  
$$
Hence the generic fiber is $\QQ(a)$, the special fiber  for $a=1$  is $\QQ[z]/\langle z^2\rangle$,
and the special fibers for $a \ne 1$ are $\QQ$. Moreover, note that the re-embeddings
$\Phi_L$ and~$\Phi_\Gamma$ are all optimal.
\end{example}

Next we turn to the question when the generic fiber or a special fiber
are isomorphic to polynomial rings.
The question to decide if an affine variety is isomorphic to an affine space is
quite delicate and connected to several hard and unsolved problems in affine geometry.
Since we study it here using the language of algebra, we introduce the following terminology.
Recall that the free objects in the category of commutative $K$-algebras are the 
polynomial rings over~$K$.

\begin{definition}\label{def-affineness}
Let $P/I$ be an affine $K$-algebra, where $P$ is a polynomial ring over~$K$
and~$I$ is an ideal in~$P$. 
\begin{enumerate}
\item[(a)] The ring $P/I$ is called a  {\bf free commutative $K$-algebra}
if there exist a polynomial ring~$P'$ over~$K$ and a $K$-algebra isomorphism 
$P/I \cong P'$. Since we are exclusively 
dealing with commutative rings here, we also call $P/I$ a {\bf free $K$-algebra}
in this case.

\item[(b)] Let $Z$ be a subtuple of the tuple of indeterminates~$X$ of~$P$.
The ring~$P/I$ is called a {\bf $Z$-separating free $K$-algebra} if there exists 
a $Z$-separating re-embedding $\Phi:\; P/I \longrightarrow K[X\setminus Z]$.  
\end{enumerate}
\end{definition}

\begin{remark}\label{rem-notenough}
If the ideal~$I$ contains a coherently $Z$-separating tuple for which the corresponding
re-embedding is optimal, and if the local ring $(P/I)_{\M/I}$ is regular, then
the affine space criterion (cf.~\cite{KLR2}, Proposition~6.7) says that $P/I$ is a $Z$-separating
free $K$-algebra. However, if the ideal~$I$ does not contain a long enough
$Z$-separating tuple, then the regularity of $(P/I)_{\M/I}$, and even the regularity of~$P/I$,
are obviously necessary conditions for $P/I$ to be a free $K$-algebra, but
they are not sufficient, as the following examples show.
\end{remark}

\begin{example}\label{ex-crachiola}
Let $P=\QQ[x,y,z,w]$ and $\M= \langle x,y,z,w\rangle$.
\begin{enumerate}
\item[(a)] The ideal $I=\langle x +x^2y +z^2 + w^3 \rangle \subseteq P$ defines 
the  {\it Russel cubic} (see~\cite{Cr}). 
There is no indeterminate for which~$I$ contains a separating polynomial.
Since $\dim_{\QQ}(\Lin_\M(I))=1$, it follows that there is no optimal separating re-embedding of~$P/I$. 
Moreover, the ring $P/I$ is regular. However, in~\cite{Cr} a proof is given, which uses 
techniques completely different from ours, to show that $P/I$ is not a free $K$-algebra.

\item[(b)] Now consider the similar looking ideal $J=\langle x + xy^2 + yz + w^3\rangle$.
Again, the ring $P/J$ is regular. Notice that~$J$ is $W\!$-homogeneous for
the non-negative grading on~$P$ defined by $W = (3,0,3,1)$.
In the next section (see Example~\ref{ex-similCrachiolaCont}) we will show how to construct
a $\QQ$-algebra isomorphism  $\Phi:\; P/I \longrightarrow \QQ[y,z,w]$
which is induced by $x \mapsto -yz-w^3$, $y\mapsto y$, $z\mapsto z + y^2z + yw^3$, 
and $w\mapsto w$. Hence the $K$-algebra $P/J$ is free, but it is not $Z$-separating free for 
any tuple of indeterminates~$Z$.

\end{enumerate}
\end{example}

In addition to the notation introduced at the beginning of this section,
we need one more abbreviation. Namely, for $\Gamma = (c_1, \dots c_m) \in K^m$, 
we let~$\M_\Gamma$ be the maximal ideal of~$P$ given by
$\M_\Gamma =\langle a_1-c_1, \dots, a_m-c_m,\, x_1, \dots, x_{n-m} \rangle$.

Our next result says that if we are in the
non-negatively graded setting and the local ring $(P/I)_{\M_\Gamma/I}$ is regular, 
then this special fiber of the family is isomorphic to an affine space.

\begin{theorem}{\bf (Freeness of a Special Fiber)}\label{thm-freefiber}\\
Let~$I$ be a $W\!$-homogeneous ideal in~$P$ such that $I\cap P_0 =\{0\}$, 
and let $\Gamma = (c_1,\dots,c_m) \in K^m$. Then the following claims hold.
\begin{enumerate}

\item[(a)] There is an isomorphism of $K$-vector spaces
$$
\alpha_\Gamma:\; \Lin_{\M_\Gamma}(I) \cong \Lin_{\langle X^\plus \rangle}(I_\Gamma)
$$
induced by the $K$-algebra homomorphism $\beta:\; P \longrightarrow K[\Xp]$ defined by
$\beta(a_i)=c_i$ for $i=1,\dots, m$ and $\beta(x_j)=x_j$ for $j=1,\dots, n-m$.

\item[(b)] We have $\dim_K ( \Cot_{\langle X^\plus \rangle}) (K[\Xp]/I_\Gamma) ) =
\dim_K (\Cot_{\M_\Gamma}(P/I) ) - m$.

\item[(c)]  If the localization $(P/I)_{\M_\Gamma/I}$ is a regular local ring, then the localization 
$(K[\Xp]/I_\Gamma)_{\langle X^\plus \rangle/I_\Gamma}$ of the fiber is a regular local ring, as well. 

In this case we have 
$\dim(K[\Xp]/I_\Gamma) = \dim\big((P/I)_{\M_\Gamma/I}\big) -m$.

\item[(d)] If the localization $(K[\Xp]/I_\Gamma)_{\langle X^\plus \rangle/I_\Gamma}$ of the fiber is 
a regular local ring, then the fiber $K[\Xp]/I_\Gamma$  is a free $K$-algebra.
\end{enumerate}
\end{theorem}

\begin{proof} 
To prove (a), we write $I = \langle g_1, \dots, g_r\rangle$, where 
$g_1,\dots, g_r\in P$ are $W\!$-homo\-geneous polynomials of positive degree. 

Every polynomial $f \in I$ can be written in the form
$f= \sum_{i=1}^{n-m} h_i\, x_i + q$, where $h_i\in K[a_1,\dots,a_m]$ and $q\in \langle \Xp\rangle ^2$. 
Consequently, the $\M_\Gamma$-linear part of~$f$ is
$$
\Lin_{\M_\Gamma}(f) = \Lin_{\M_\Gamma}( \tsum_{i=1}^{n-m} h_i \, x_i) = \tsum_{i=1}^{n-m} h_i(c_1,\dots,c_m)\, x_i
$$
To obtain the fiber ideal~$I_\Gamma$, we have to substitute $a_i \mapsto c_i$ for $i=1,\dots,m$
in the polynomials in~$I$. The $\langle X^+\rangle$-linear parts of the resulting polynomials are precisely the polynomials
$\Lin_{\M_\Gamma}(f)$, viewed as polynomials in $K[\Xp]$. This observation implies claim~(a).

Now let us prove~(b).  We have 
\begin{align*}
\dim_K ( \Cot_{\langle X^\plus \rangle} (K[\Xp]/I_\Gamma) ) &\;=\; 
n - m - \dim_K (\Lin_{\langle X^\plus \rangle}(I_\Gamma)) \cr
&\; = \; n - m - \dim_K ( \Lin_{\M_\Gamma}(I) )\cr 
&\;=\; \dim_K ( \Cot_{\M_\Gamma}(P/I)  ) - m
\end{align*}
where the second equality follows from~(a).

Next we show~(c). Let $S=(K[\Xp]/I_\Gamma)_{\langle X^\plus \rangle/I_\Gamma}$.
We claim that there is the following chain of relations.
\begin{align*}
\dim((P/I)_{\M_\Gamma/I} ) - m  &\;\le\;  \dim(S) \;\le\; 
\dim_K ( \Cot_{\langle X^\plus \rangle/I_\Gamma}(K[\Xp]/I_\Gamma) ) \cr
&\;=\; \dim_K ( \Cot_{\M_\Gamma}(P/I) ) - m  \;=\; \dim((P/I)_{\M_\Gamma/I} ) - m
\end{align*}
The first inequality follows from the fact that~$S$ is obtained from $(P/I)_{\M_\Gamma/I}$ 
by modding out the ideal~$\langle a_1-c_1, \dots, a_m-c_m \rangle$ which is generated by~$m$ elements. 
The second inequality follows from the fact that the dimension of a local ring is bounded above by the 
dimension of the  cotangent space at its maximal ideal (see~\cite{Kun}, V.4.5). 
The first  equality follows from part~(b), and the second equality is a consequence of the assumption that  
$(P/I)_{\M_\Gamma/I}$ is regular.

It follows that all relations are in fact equalities. In particular, we have 
$\dim(S) = \dim_K (\Cot_{\langle X^\plus \rangle/I_\Gamma} (K[\Xp]/I_\Gamma))$, and hence~$S$ 
is a regular local ring. Notice that we have $\dim(K[\Xp]/I_\Gamma) = \dim((K[\Xp]/I_\Gamma)_{\langle \Xp\rangle
/ I_\Gamma})$ for the positively graded ring $K[\Xp]/I_\Gamma$.

To prove~(d), we use~(c) and Proposition~\ref{prop-OptFibers} to deduce that 
the assumptions of the affine space criterion (cf.~\cite{KLR2}, Proposition~6.7) are satisfied.
As already pointed out in Remark~\ref{rem-notenough}, the claim follows.
\end{proof}

Notice that the analogous statement for the generic fiber, namely that $L[\Xp]/I_L$ is a free $L$-algebra
if it is a regular ring, follows from the affine space criterion (cf.~\cite{KLR2}, Proposition~6.7), as well.

The proof of part~(b) of this theorem shows that the images of 
$a_1-c_1,\dots, a_m-c_m$ 
in the regular local ring $(P/I)_{\M_\Gamma/I}$ are part of a
regular system of parameters. Thus~\cite{Kun}, VI.1.10, is another way to deduce claim~(c).

\medskip
Finally, we note that the scheme $\Spec(P/I)$ is connected in our setting.

\begin{proposition}\label{prop-isprime}
Let~$I$ be a $W\!$-homogeneous ideal in~$P$ such that $I\cap P_0 =\{0\}$.
\begin{enumerate}
\item[(a)] The scheme $\Spec(P/I)$ is connected.

\item[(b)] If $P/I$ is a regular ring, then~$I$ is a prime ideal.

\item[(c)] If $P/I$ is a regular ring, then for every $\Gamma\in K^m$ the fiber
$K[\Xp]/I_\Gamma$ has dimension $\dim(P/I)-m$.

\end{enumerate}
\end{proposition}

\begin{proof}
To show~(a), we note that $P_0$ is a polynomial ring and the special fibers of the 
surjective morphism $\Spec(P/I) \longrightarrow \Spec(P_0)$ are connected, because
the rings $K[\Xp] / I_\Gamma$ are positively graded.

Claim~(b) follows from~(a) and the well-known fact that a regular local ring
is an integral domain. 

Lastly, we show~(c). By~(b), the ring $P/I$ is equidimensional. Hence we have
$\dim(P/I) = \dim((P/I)_{\M_\Gamma/I})$ and the claim follows from part~(c) of the
preceding theorem.
\end{proof}

In the language of Algebraic Geometry, Theorem~\ref{thm-freefiber} and this proposition
imply that the family $\Spec(P/I) \longrightarrow \Spec(P_0)$ is an $\mathbb{A}^k$-fibration
for $k=\dim(P/I)-\dim(P_0)$ (cf.~\cite{Gup}, Section~4).

\bigskip\bigbreak
%
%

\section{Re-Embeddings Using the Unimodular Matrix Problem}
\label{sec6}

In this section we introduce a further technique for constructing re-embeddings.
It is based on the unimodular matrix problem which is defined as follows.

\begin{definition}\label{def-UMP}
Let $P$ be a polynomial ring over a field~$K$, and let $1\le k < \ell$.
\begin{enumerate}
\item[(a)] A matrix $A \in \Mat_{k,\ell}(P)$ is said to be {\bf unimodular} 
if its maximal minors generate the unit ideal.
 
\item[(b)] A matrix $B\in \Mat_\ell(P)$ is said to solve the {\bf unimodular matrix problem (UMP)} for~$A$
if $\det(B)=1$ and $A\cdot B = (I_k \mid 0)$, where $I_k$ denotes the identity matrix of size~$k$.

\end{enumerate}
\end{definition}

This definition can clearly be extended to the case $k=\ell$. 
In this case we simply have $B = A^{-1}$.
The following remark provides another view of this problem.

\begin{remark}\label{rem-UMC}
Given~$B$ which solves the UMP for~$A$, 
we have $A = (I_k \mid 0) \cdot B^{-1}$. In other words, the first~$k$ rows of
the matrices~$A$ and $B^{-1}$ coincide. For this reason the matrix $B^{-1}$ is
called a {\bf unimodular matrix completion} of~$A$.
\end{remark}

The following  example shows that the solution of a UMP is not unique.

\begin{example}\label{ex-UMC}
Let $P = \QQ[x]$, and let $A = \begin{pmatrix} 1 +x^2  & x \end{pmatrix}$.
Then, among others, we get the UMP solutions
$B = \left( \begin{smallmatrix} 1 & -x\;  \\ -x \; & 1+x^2 \end{smallmatrix} \right)$
and $\tilde{B} = \left( \begin{smallmatrix} 1-x^2 & -x \\ x^3 & 1+x^2 \end{smallmatrix} \right)$.
\end{example}

In general, the computation of a UMP solution is a difficult task, but effective algorithms
are available and have been described for instance in~\cite{F}, \cite{FG}, and~\cite{LS}. 
Many tricks how to shortcut the computation in special cases are contained in~\cite{FQ}.
A description of an implementation in \texttt{Macaulay2} is presented in~\cite{BS}.
Based on the calculation of UMP solutions, one has effective methods for computing 
a basis of a finitely generated projective $P$-module given as the kernel or the cokernel of a
polynomial matrix (see~\cite{LS}, Section~3), thereby making the famous Quillen-Suslin theorem
(cf.~\cite{Kun}, IV.3.15) effectively computable.
For our purposes, we proceed as follows.

\begin{algorithm}{\bf (Computing a UMP Solution)}\label{alg-completion}\\
Let $P$ be a polynomial ring over a field~$K$,
let $k<\ell$, and let $A=(h_{ij})$ be a  unimodular matrix  in $\Mat_{k,\ell}(P)$.
For $j=1,\dots, \ell$, let $h_j^\ast \in P^k$ be the $j$-th column of~$A$.
Consider the following instructions:
\begin{enumerate}
\item[(1)] For $j=1,\dots, k$, use explicit module membership to calculate
polynomials $c_{1j},\dots,c_{\ell j} \in P$ such that $e_j = c_{1j} h_1^\ast + \cdots
+ c_{\ell j} h_\ell^\ast$.

\item[(2)] Form the matrix $C = (c_{ij}) \in \Mat_{\ell,k}(P)$.

\item[(3)] Compute a basis $\{b_1,\dots,b_{\ell-k}\} \subseteq P^\ell$ of the
syzygy module $\Syz_P(h_1^\ast,\dots, h_\ell^\ast)$.

\item[(4)] Append the columns $b_1,\dots,b_{\ell-k}$ to the matrix~$C$ and obtain
a matrix $B' \in \Mat_\ell(P)$. Then divide the last column of~$B'$ 
by $\det(B')$ to get a matrix $B\in \Mat_\ell(P)$ with $\det(B)=1$.

\item[(5)] Return the matrix~$B$ and stop.
\end{enumerate}
This is an algorithm which computes a matrix $B\in \Mat_\ell(P)$
which solves the UMP for~$A$.
\end{algorithm}

\begin{proof}
Step~(1) can be executed, since the $P$-linear map $\phi:\; P^\ell \longrightarrow
P^k$ defined by~$A$ is surjective, as the $P$-module $P^k/U$, where~$U$ is 
generated by the columns of the matrix~$A$, has its 0-th Fitting ideal 
$F_0(P^k/U) = \langle 1\rangle$ and is hence the zero module (cf.~\cite{Eis}, Prop.~20.8).

Step~(3) can be executed, because the module $\Ker(\phi)$ is free by the 
Quillen-Suslin theorem (see \cite{Kun}, IV.3.15). 
Step~(4) can be executed, because the columns of~$B'$ are a $P$-basis of~$P^\ell$,
since the exact sequence
$$
0 \;\longrightarrow\; \Syz_P(h_1^\ast,\dots,h_\ell^\ast) \;\longrightarrow\; P^\ell  \;\longrightarrow\;
P^k \;\longrightarrow\; 0
$$
splits. Thus we have $\det(B')\in K \setminus \{0\}$.
Hence the algorithm is finite and executable, and it remains to show correctness.

By construction, we have $A\cdot C = I_k$, where $I_k$ is the identity matrix of size~$k$.
Then we obtain $A\cdot B = (I_k \mid 0) \in \Mat_{k,\ell}$, and thus~$B$ solves the UMP for~$A$.
\end{proof}

The hard part here is step~(3) which may require us to go through the full
algorithm for the Quillen-Suslin theorem. 

\medskip
Now we are ready for the first main result of this section
which provides us with a new method for constructing re-embeddings in the case 
when a certain Fitting ideal of some coefficient matrix is the unit ideal.
So, we return to the setting of the preceding section, in particular 
we let $P=K[a_1,\dots,a_m,x_1,\dots,x_{n-m}]$
be non-negatively graded by $W=(w_1,\dots,w_n)$, assume that $w_1\le \cdots \le w_n$,  that
$P_0=K[a_1,\dots,a_m]$, and let $\Xp = (x_1,\dots,x_{n-m})$.
In this setting, we use the solution of the UMP to define a new type of re-embedding as follows.

\begin{algorithm}{\bf (Re-Embeddings Using the UMP)}\label{alg-UMC-reemb}\\
Let $g_1,\dots,g_r \in P$ be homogeneous polynomials of positive $W\!$-degree, 
and let $I=\langle g_1,\dots,g_r\rangle$.
Assume that there exists a number $k$ with $1\le k\le \min\{r,n-m\}$
such that, if we write $g_i = h_{i1} x_1 + \cdots + h_{i\, n-m} x_{n-m} + \tilde{h}_i$ with $h_{ij} \in P_0$
and $\tilde{h}_i \in \langle\Xp\rangle^2$ for $i=1,\dots,k$, then the maximal minors of the matrix 
$A = (h_{ij}) \in \Mat_{k,n-m}(P_0)$  generate the unit ideal. Moreover, let $g_1,\dots,g_k$ be ordered such that
$\deg_W(g_1) \le \cdots \le \deg_W(g_k)$, and let $D=\{\deg_W(x_i) \mid i=1,\dots,n-m\}$.
Consider the following sequence of instructions.
\begin{enumerate}
\item[(1)] For every $d\in D$, in increasing order, perform the following steps (2)-(5).

\item[(2)] Let $1\le \alpha\le\beta \le n-m$ be such that $x_\alpha, x_{\alpha+1},\dots,x_\beta$
are the indeterminates of degree~$d$, and let $1\le \gamma \le \delta \le k$ be such that
$g_\gamma,\dots,g_\delta$ are the polynomials of degree~$d$ among $\{g_1,\dots,g_k\}$.

\item[(3)] Let $A_d = (h_{ij}) \in \Mat_{\delta -\gamma +1,\, \beta-\alpha+1}(P_0)$ be the submatrix
corresponding to the block of~$A$ in degree~$d$.

\item[(4)] Perform Algorithm~\ref{alg-completion} and compute a matrix $B_d \in \Mat_{\beta-\alpha+1}(P_0)$
which solves the UMP for~$A_d$.

\item[(5)] Let $F_d = (x_\alpha,\dots,x_\beta) \cdot B_d^{\rm tr}$.

\item[(6)] Let $F$ be the concatenation of the tuples $(a_1,\dots,a_m)$
and $F_d$ with $d\in D$, and let~$\phi$ be the 
$K\!$-algebra endomorphism of~$P$ defined by~$F$.

\item[(7)] For $i=1,\dots,k$, write $\phi(g_i) = x_ {\nu_i} - q_i$ with $q_i\in P$ 
such that~$x_{\nu_i}$ divides no term in $\Supp(q_i)$. Rewrite~$\phi(g_i)$ using 
$(x_{\nu_1} - q_1, \dots, x_{\nu_{i-1}} - q_{i-1})$ 
to obtain a polynomial $x_{\nu_i}-\tilde{q}_i$.

\item[(8)] Let $\widehat{X} = \{x_1,\dots,x_{n-m}\} \setminus \{x_{\nu_1},\dots,x_{\nu_k}\}$,
and let $\widehat{P} = P_0[\widehat{X}]$ be graded by the restriction $\widehat{W}$ of~$W\!$.
Define a tuple $\widehat{F}=(a_1,\dots,a_m, f_1,\dots, f_{n-m})$ of polynomials in~$\widehat{P}$
where $f_{\nu_i} = \tilde{q}_i$ for $i=1,\dots,k$ and $f_j=x_j$ for $x_j\in\widehat{X}$.
Let $\psi:\; P \longrightarrow \widehat{P}$ be the $K\!$-algebra homomorphism given by~$\widehat{F}$.

\item[(9)] To define the composition $\theta = \psi\circ\phi:\; P \longrightarrow \widehat{P}$, substitute
the entries of~$\widehat{F}$ for the indeterminates in~$F$ and return the resulting tuple~$T$.
\end{enumerate}
This is an algorithm which computes a tuple of polynomials~$T$ in~$\widehat{P}$ such that
the corresponding $K\!$-algebra homomorphism $\theta:\; P \longrightarrow \widehat{P}$
induces a $K\!$-algebra isomorphism $\Theta:\; P/I \longrightarrow \widehat{P}/J$.
Here $\widehat{P} = P_0[\widehat{X}]$ is a polynomial ring in $n-k$ indeterminates 
and $J = \langle \theta(g_{k+1}), \dots, \theta(g_r)\rangle$ is a $\widehat{W}$-homogeneous ideal in~$\widehat{P}$.
\end{algorithm}

\begin{proof}
First of all, notice that the homogeneity of the polynomials $g_1,\dots,g_k$ and their ordering by
non-decreasing $W\!$-degrees implies that the matrix~$A$ is block diagonal, where each block~$A_d$ has
at least as many columns as rows and where the maximal minors of~$A_d$ generate the unit ideal in~$P_0$. 

For every $d\in D$, Algorithm~\ref{alg-completion} shows that the matrix~$B_d$ satisfies
$$
A_d \cdot B_d \cdot (x_\alpha,\dots,x_\beta)^{\rm tr} \;=\; (I_{\delta-\gamma+1} \mid 0) \cdot
(x_\alpha,\dots,x_\beta)^{\rm tr}  \;=\; \begin{pmatrix} x_{\nu_\gamma} && 0 & 0 \cdots 0 \\
& \ddots & & \vdots \qquad \vdots \\
0 && x_{\nu_\delta} & 0 \cdots 0  \end{pmatrix}
$$
where $\nu_\gamma = \alpha$, $\dots$, $\nu_\delta = \alpha + \delta-\gamma$.
Since we have $F_d = (x_\alpha,\dots,x_\beta)\cdot B_d^{\rm tr} = (\phi(x_\alpha),\dots, \phi(x_\beta))$, we obtain
$$
\phi(g_i) \;=\; \phi(h_{i\alpha} x_\alpha + \cdots + h_{i \beta} x_\beta) + \phi(\tilde{h}_i) \; = 
x_{\nu_i} + \phi(\tilde{h}_i)
$$
for $i=\gamma,\dots, \delta$. Notice that this implies that~$\phi$ is homogeneous of degree zero.

Thus the polynomials $\phi(g_1),\dots\phi(g_k)$ have the shape required by step~(7) and
are $W$-homogeneous of non-decreasing degree. As $\tilde{h}_i\in \langle \Xp\rangle^2$, we also get that
the tuple $(\phi(g_1),\dots,\phi(g_k))$ is $Z$-separating for $Z=(x_{\nu_1}, \dots, x_{\nu_k})$.
Therefore step~(7) is nothing but the conversion of this tuple to a coherently $Z$-separating tuple
(see Algorithm~\ref{alg-findCohF}),
and the map~$\psi$ induces the $Z$-separating re-embedding of the ideal $\phi(I)$.

As the matrix~$B$ is invertible, it follows that~$\phi$ is a $K\!$-algebra automorphism of~$P$.
Consequently, the map~$\Theta$ is a $K\!$-algebra isomorphism.
Since both~$\phi$ and~$\psi$ are homogeneous of degree zero, the ideal~$J$ is $\widehat{W}$-homogeneous.
\end{proof}

Let us apply this algorithm to a concrete example.

\begin{example}\label{ex-twodeg}
Let $P =\QQ[a, x_1, x_2, x_3, x_4, x_5, x_6 ]$ be
non-negatively graded by $W=(0, 1, 2, 2, 2, 4, 4)$,  let 
$G = (g_1, g_2, g_3)$, where  
\begin{align*}
g_1 &\;=\; ax_2 +(a^2-1)x_3 + x_1^2 \\
g_2 &\;=\; (a-1)x_5 +a^2x_6 + a^3 x_1^2x_2 +x_1^2x_4 \\
g_3 &\;=\;  x_2x_4 -x_3^2 + a x_5
\end{align*}
are $W\!$-homogeneous polynomials, and let $I=\langle g_1,g_2,g_3\rangle$.

The polynomials $g_1$ and $g_2$ have representations such that the maximal minors of
the matrix 
$$
A \;=\; \begin{pmatrix}
0 & a & a^2-1 & 0 & 0 & 0\\
0 & 0 & 0 & 0 & a-1 & a^2
\end{pmatrix}
$$
generate the unit ideal. Hence the hypotheses of Algorithm~\ref{alg-UMC-reemb}
are satisfied and we can follow its steps.
\begin{enumerate}
\item[(1)] We have $D=\{1,2,4\}$.

\item[(2)-(5)] For $d=1$, the matrix $A_1$ is empty and there is nothing to do.

\item[(2)-(5)] For $d=2$, we have $A_2 = ( a\;\; a^2{-}1 \;\;  0)$. This yields
$$
B_2 = \begin{pmatrix} 
a & a^2-1 & 0\\
-1 & -a & 0\\
0 & 0 & -1
\end{pmatrix}
$$
We check that $\det(B_2) =1$ and deduce that the tuple~$F_2$ is given by
$F_2=(a x_2 + (a^2-1) x_3,\; -x_2 -ax_3,\; -x_4)$.

\item[(2)-(5)] For $d=4$, we have $A_4 = (a{-}1 \;\; a^2)$. This yields
$B_4 = \left( \begin{smallmatrix} -a{-}1 \ & -a^2\\ 1\  & a{-}1 \end{smallmatrix} \right)$.
We check that $\det(B_4) =1$ and deduce that the tuple~$F_4$ is given by
$F_4 = ( (-a{-}1) x_5 -a^2 x_6,\; x_5 + (a -1)x_6 )$.

\item[(6)] Altogether, we get the concatenated tuple
$$
F = (a,\, x_1,\, a x_2 +(a^2-1) x_3,\, -x_2 -ax_3,\, -x_4,\, (-a-1) x_5 - a^2 x_6,\,
x_5 + (a-1) x_6)
$$
which defines a $\QQ$-algebra automorphism~$\phi$ of~$P$.

\item[(7)] Now we apply~$\phi$ to $g_1,g_2$ and we get
$\phi(g_1) = x_2 + x_1^2$ and $\phi(g_2) = x_5 + a^4 x_1^2 x_2 + (a^5-a^3) x_1^2 x_3 - x_1^2 x_4$.
Hence we have $q_1 = \tilde{q}_1 = -x_1^2$, 
$q_2 = -a^4 x_1^2 x_2 - (a^5-a^3) x_1^2 x_3 + x_1^2 x_4$,
and rewriting $x_2 \mapsto -x_1^2$ yields $\tilde{q}_2 = a^4 x_1^4 - (a^5-a^3) x_3 + x_1^2 x_4$.

\item[(8)] Thus we get $\widehat{F} = (a,\, x_1,\, -x_1^2,\, x_3,\, x_4,\, a^4 x_1^4 - (a^5-a^3) x_3 + x_1^2 x_4,\, x_6)$.

\item[(9)] Next we substitute the entries of~$\widehat{F}$ for the indeterminates
in the tuple~$F$ and get a tuple~$T$ which defines a $K\!$-algebra homomorphism $\theta:\;
P \longrightarrow \widehat{P}$. We return~$T$ and stop.
\end{enumerate}

The homomorphism~$\theta$ induces a $K\!$-algebra isomorphism
$\Theta:\; P/I \longrightarrow \widehat{P}/J$ where $\widehat{P} = 
\QQ[a,x_1,x_3,x_4,x_6]$ and $J= \langle \theta(g_3)\rangle$ with
\begin{gather*}
\theta(g_3) \;=\; a^7x_1^2x_3 -a^6x_1^4 +a^6x_1^2x_3 -a^5x_1^4 -a^5x_1^2x_3 -a^4x_1^2x_3  -a^3x_6  \\
  -a^2x_3^2 -a^2x_1^2x_4 -a^2x_3x_4 +2ax_1^2x_3 -x_1^4 +x_3x_4
\end{gather*}
Since the hypersurface defined by~$\theta(g_3)$ is singular at the origin, this is an optimal re-embedding of~$I$.

\end{example}

Finally, we examine when the re-embeddings constructed using the
preceding algorithm allow us to conclude that the given $K$-algebra is free.
The next theorem is the key result.

\begin{theorem}{\bf (Freeness of Regular Non-Negatively Graded Algebras)}\label{thm-Free}\\
Let $K$ be a perfect field, let $P=K[a_1,\dots,a_m,x_1,\dots,x_{n-m}]$ be non-negatively graded by~$W\!$, 
where $P_0=K[a_1,\dots,a_m]$,
and let~$I$ be a $W\!$-homogeneous ideal in~$P$ with $I\cap P_0 = \{0\}$.
If $P/I$ is a regular ring, then the following claims hold.
\begin{enumerate}
\item[(a)] The ideal~$I$ is generated by a $W\!$-homogeneous regular sequence
$f_1,\dots,f_k$, where $k=n - \dim(P/I)$.

\item[(b)] When we write $f_i = p_{i1} x_1 + \cdots + p_{i n-m}x_{n-m} + q_i$
with $p_{ij}\in P_0$ and $q_i \in \langle \Xp\rangle^2$, then the maximal minors of the matrix
$A = (p_{ij}) \in \Mat_{k,n-m}(P_0)$ generate the unit ideal.

\item[(c)] There exist a polynomial ring $\widehat{P} = P_0[\widehat{X}]$ in $n-k$ indeterminates
and a homogeneous $K$-algebra isomorphism $P/I \cong \widehat{P}$. In particular, the ring 
$P/I$ is a free $K$-algebra.
\end{enumerate}
\end{theorem}

\begin{proof}
First we show~(a) and~(b). Let $\{g_1,\dots,g_r\}$ be a $W\!$-homogeneous system 
of generators of~$I$, and assume that $\deg_W(g_1) \le \cdots \le \deg_W(g_r)$.
For $i=1,\dots,r$, we write $g_i = h_{i1} x_1 + \cdots + h_{i n-m} x_{n-m} + \tilde{q}_i$
with $h_{ij}\in P_0$ and $\tilde{q}_i\in \langle \Xp\rangle^2$.
Moreover, let $H=(h_{ij}) \in \Mat_{r,n-m}(P_0)$, and let $k = n - \dim(P/I)$.

As an {\it intermediate claim}, we prove that the minors of size $k+1$ of~$H$ are zero
and the ideal~$J$ generated by the minors of size~$k$ of~$H$ is the unit ideal.
Let $\overline{K}$ be the algebraic closure of~$K$ and
$\overline{P} = \overline{K}[a_1,\dots,a_m,x_1,\dots,x_{n-m}]$.
Since~$K$ is perfect and $P/I$ is regular, also the ring $\overline{P} / I\, \overline{P}$
is regular. Furthermore, by Proposition~\ref{prop-isprime}, we know that~$I$ is a prime ideal,
and hence also the ideal $I\,\overline{P}$ is a prime ideal. Additionally, we have
$k = n - \dim(\overline{P}/I\, \overline{P})$.

Let $\Gamma = (c_1,\dots,c_m) \in \overline{K}^{\,m}$ and 
$\M_\Gamma = \langle a_1-c_1,\dots,a_m-c_m,x_1,\dots,x_{n-m}\rangle$ in~$\overline{P}$.
Then the local ring $(\overline{P}/I\,\overline{P})_{\M_\Gamma / I\,\overline{P}}$ is regular,
and the Jacobian criterion (see~\cite{Kun}, VI.1.5) implies that
the rank of the Jacobian matrix at~$\Gamma$  is 
$k_\Gamma =  n - \dim((\overline{P}/I\,\overline{P})_{\M_\Gamma / I\,\overline{P}})$.
As $\overline{P}/I\,\overline{P}$ is a regular integral domain, we have $k_\Gamma=k$ here.

Moreover, since  $(\partial g_i/\partial x_j)(c_1,\dots,c_m,0,\dots,0) = h_{ij}(c_1,\dots,c_m)$
for $i = 1, \dots, r$ and $j =1, \dots, n-m$,
and since $(\partial g_i/\partial a_i)(c_1,\dots,c_m,0,\dots,0) = 0$ for $i = 1,\dots,r$ and  $j=1,\dots, n-m$,
we deduce that  the rank of the Jacobian matrix evaluated at the point $\Gamma$ 
is also the rank of the matrix $H_\Gamma$, i.e., the rank of the matrix~$H$ evaluated at~$\Gamma$.
Hence the rank of $H_\Gamma$ is~$k$ for every point $\Gamma  \in \overline{K}^{\,m}$.

Consequently, the minors of size $k+1$ of~$H$ evaluate to zero at every point~$\Gamma$,
i.e., they are zero, and the ideal~$J$ generated by the minors of size~$k$ of~$H$
evaluates to a non-zero ideal $J_\Gamma$ at every point $\Gamma  \in \overline{K}^{\,m}$.
Thus the ideal~$J$ is the unit ideal in $\overline{K}[a_1,\dots,a_m]$.
Since it is generated by polynomials defined over~$K$, it follows that~$J$  
is the unit ideal in $K[a_1,\dots,a_m]$ and the {\it intermediate claim} is proved.

Using~\cite{Eis}, Proposition~20.8, we now conclude that the $P_0$-submodule~$U$ of~$P_0^{n-m}$ 
generated by the rows of~$H$ is locally free of rank~$k$. Then the Quillen-Suslin theorem
implies that~$U$ is a free $P_0$-module of rank~$k$. Let $\{v_1,\dots,v_k\}$ be a $P_0$-basis
of~$U$, and write $v_i = \sum_{j=1}^r u_{ij}\cdot (h_{j1},\dots,h_{j\, n-m})$ with $u_{ij}\in P_0$
for $i=1,\dots,k$. 

Next we consider the polynomials $f_i = \sum_{j=1}^r u_{ij}\, g_j \in I$.
Notice the the entries $h_{ij}$ of~$H$ can only be non-zero if $\deg_W(g_i) = \deg_W(x_j)$.
W.l.o.g.\ we order the indeterminates such that $\deg_W(x_1) \le \cdots \le \deg_W(x_{n-m})$.
Then the ordering of the generators~$g_i$ and of the indeterminates~$x_j$ by 
non-decreasing $W\!$-degree yields a block diagonal structure of~$H$. 
Hence the polynomials~$f_i$ will be $W\!$-homogeneous, too, if the 
basis $\{v_1,\dots,v_k\}$ is chosen accordingly.

Now the ideal $\langle f_1,\dots,f_k\rangle$ is $W\!$-homogeneous and contained in~$I$.
When we write $f_i = p_{i1} x_1 + \cdots + p_{i\, n-m}x_{n-m} + q_i$
with $p_{ij}\in P_0$ and $q_i \in \langle \Xp\rangle^2$ as in~(b), then the vectors
$v_i = (p_{i1},\dots,p_{i\, n-m})$ for $i=1,\dots,k$ are a $P_0$-basis of the row
space of the matrix $(p_{ij}) \in \Mat_{k,n-m}(P_0)$.
Consequently, the maximal minors of the matrix $(p_{ij})$ generate the unit ideal
in~$P_0$, and the Jacobian criterion implies that $P/\langle f_1,\dots,f_k\rangle$
is a smooth homogeneous complete intersection of dimension $n-k = \dim(P/I)$.
Since the local rings of this complete intersection are integral domains (cf.~\cite{Kun}, V.5.15.a),
the corresponding local rings of~$P/I$, which are residue class rings of the same
dimension, have to agree with them. By the local-global principle (cf.~\cite{Kun}, IV.1.4), 
it follows that $I=\langle f_1,\dots,f_k\rangle$.

To show~(c), we now apply Algorithm~\ref{alg-UMC-reemb}. Its hypothesis is 
satisfied by~(b). The result is a homogeneous $K$-algebra isomorphism
$\Theta:\; P/I \longrightarrow \widehat{P}/J$, where $\widehat{P}=P_0[\widehat{X}]$
has dimension $\dim(\widehat{P})= n-k = \dim(P/I)$ and where~$J$ is a homogeneous 
ideal in~$\widehat{P}$ which is necessarily zero.
\end{proof}

Using Algorithm~\ref{alg-UMC-reemb}, we can turn the proof of this theorem 
into the following effective method for computing an isomorphism between $P/I$ and a polynomial ring.
To perform step~(2) of the following algorithm, one can use the method given in~\cite{LS}, Sect.~3
to compute a basis of a free submodule of a free module over a polynomial ring.

\begin{algorithm}{\bf (Re-Embedding a Regular Graded Algebra)}\label{alg-free}\\
In the setting of the theorem, let $\{g_1,\dots,g_r\}$ be a $W\!$-homogeneous system of generators
of~$I$ such that $\deg_W(g_1) \le \cdots \le \deg_W(g_r)$.
Consider the following sequence of instructions.
\begin{enumerate}
\item[(1)] For $i=1,\dots,r$, write $g_i = h_{i1} x_1 + \cdots + h_{i\, n-m} x_{n-m} +
\tilde{h}_i$ with $h_{ij} \in P_0$ and $\tilde{h}_i \in \langle \Xp \rangle^2$. 
Form the matrix $H=(h_{ij}) \in \Mat_{r,n-m}(P_0)$.

\item[(2)] Compute a $P_0$-basis
$\{v_1,\dots,v_k\}$ of the row space of~$H$. For $i=1,\dots,k$, write $v_i= \sum_{j=1}^r u_{ij}\cdot
(h_{j1},\dots,h_{j\,n-m})$ with $u_{ij}\in P_0$.

\item[(3)] Form the polynomials $f_i= \sum_{j=1}^r u_{ij} g_j$ and apply Algorithm~\ref{alg-UMC-reemb}
to compute a tuple $T = (a_1, \dots, a_m, p_1,\dots, p_{n-m})$ which defines
a homogeneous $K\!$-algebra isomorphism $\Theta:\; P/\langle f_1,\dots,f_k\rangle 
\longrightarrow \widehat{P}$, where $\widehat{P} = P_0[\widehat{X}]$ is graded by $\widehat{W}$.

\item[(4)] Return the tuples $\widehat{W}$ and $(p_1,\dots,p_{n-m})$ and stop. 
\end{enumerate}
This is an algorithm which computes tuples $\widehat{W}$ and $(p_1,\dots,p_{n-m})$ such
that there is a homogeneous $K$-algebra isomorphism 
$\Theta:\; P/I \longrightarrow \widehat{P}$ which is given as follows.

The ring $\widehat{P} = P_0[\widehat{X}]$ is non-negatively graded by~$\widehat{W}$, 
and~$\Theta$ is induced by the $K$-algebra
homomorphism $\theta:\; P \longrightarrow \widehat{P}$ such that we have 
$\theta(a_i)=a_i$ for $i=1,\dots,m$ and $\theta(x_j)=p_j$ for $j=1,\dots,n-m$.
\end{algorithm}

Let us apply this algorithm to a concrete example.
\begin{example}\label{ex-freeQalgebra}
Let $P = \QQ[a_1 , a_2, x_1, x_2, x_3, x_4, x_5]$ be non-negatively graded by $W = (0,0,1,1,1,2,2)$, let $G = (g_1,g_2)$, 
where
\begin{align*}
g_1 &= (a_1^2a_2^2 + a_1a_2 + 1)x_1 -(a_1a_2^2 + 2a_2)x_2 + (a_1a_2^3 + a_2^2)x_3\\
g_2 &= (1 -2a_1a_2) x_4 + a_1a_2 x_5 + a_1a_2x_2^2 + 2a_1x_1x_2
\end{align*}
are $W$-homogeneous polynomials, and let $I = \langle g_1, g_2\rangle$.

The polynomials $g_1$ and $g_2$ have representations such that the maximal minors
of the matrix
$$
A=\begin{pmatrix} 
0 &0 &a_1^2a_2^2 + a_1a_2 + 1 &  -(a_1a_2^2 + 2a_2)  & a_1a_2^3 + a_2^2 &0  &0\\
0 &0& 0& 0& 0& 1 -2a_1a_2&  a_1a_2         
\end{pmatrix}
$$
generate the unit ideal. As we saw in the proof of the theorem, this implies that~$P/I$ is regular.
Hence the hypotheses of Algorithm ~\ref{alg-free} are satisfied and we can follow its steps. 
As the $2\times 2$-minors of~$A$ generate the unit ideal, its rows are a basis of its row space.
Thus we immediately move to step~(3), i.e., we apply Algorithm~\ref{alg-UMC-reemb}.

(1)  We have $D = \{1,2\}$.

(2)-(5) For $d = 1$, we have 
$$
A_1= \begin{pmatrix}a_1^2a_2^2 + a_1a_2 + 1  &  -(a_1a_2^2 + 2a_2)  &a_1a_2^3 + a_2^2 \end{pmatrix}
$$
This yields
$$
B_1=\left( \begin{array}{ccc} 
  1 &  0 &  -a_{2} \\
  \frac{1}{2} a_{1}^2 a_{2} +\frac{1}{2} a_{1} & a_{1} a_{2}^2  +a_{2} & \frac{1}{2} a_{1} a_{2} -\frac{1}{2} \\
  \frac{1}{2} a_{1}^2  & a_{1} a_{2} +2 & {\frac{3}{2}}^{\mathstrut} a_{1}
\end{array}\right) 
$$
We check that $\det(B_1) = 1$ and deduce that the tuple $F_1=(f_1,f_2,f_3)$ is given by
\begin{align*}
f_1 &= x_1 -a_2x_3 \\
f_2 &=  ( \tfrac{1}{2} a_{1}^2 a_{2} +\tfrac{1}{2} a_{1} ) x_{1} +(a_{1} a_{2}^2  +a_{2} )x_2 + (\tfrac{1}{2} a_{1} a_{2} -\tfrac{1}{2})x_3\\
f_3 &= \tfrac{1}{2} a_{1}^2 x_{1} +(a_{1} a_{2} +2) x_{2} +\tfrac{3}{2} a_{1} x_{3} 
\end{align*}

(2)-(5) For $d = 2$, we have $A_2 = \begin{pmatrix} -2a_1a_2 +1 & a_1a_2 \end{pmatrix}$.
This yields
$$
B_2=\left( \begin{array}{ccc} 
  1 &  -a_{1} a_{2} \\
  2 & -2 a_{1} a_{2} +1
\end{array}\right) 
$$
We check that $\det(B_2) = 1$ and deduce that the tuple $F_2=(f_4,f_5)$ is given by
\begin{align*}
f_4 &=  x_4 -a_1a_2 x_5\\
f_5 &=  2x_4 + (-2a_1a_2+1) x_5
\end{align*}

(6) Altogether, we get the tuple $F=(a_1,a_2,f_1,f_2,f_3,f_4,f_5)$ which defines a $\QQ$-algebra
automorphism~$\phi$ of~$P$.

(7) We get $\phi(g_1) = x_1$. Then we rewrite the polynomial $\phi(g_2)$ by substituting $x_1 \mapsto 0$ and get 
$x_4- \tilde{q}$, where 
\begin{align*}
\tilde{q} = & -a_{1}^3 a_{2}^5 x_{2}^2  -a_{1}^3 a_{2}^4 x_{2} x_{3} -2 a_{1}^2 a_{2}^4 x_{2}^2  -\tfrac{1}{4} a_{1}^3 a_{2}^3 x_{3}^2
+ 2 a_{1}^2 a_{2}^3 x_{2} x_{3}\\ 
& - a_{1} a_{2}^3 x_{2}^2  +\tfrac{3}{2} a_{1}^2 a_{2}^2 x_{3}^2  +3 a_{1} a_{2}^2 x_{2} x_{3} 
-\tfrac{5}{4} a_{1} a_{2} x_{3}^2
\end{align*}

(8) Thus we get $\widehat{F} = (a_1,a_2, 0, x_2, x_3, \tilde{q}, x_5)$.

(9) Next we substitute the entries of~$\widetilde{F}$ for the indeterminates in~$F$ and get a tuple~$T$
which we return.

Altogether, we obtain an isomorphism $\Theta:\; P/I \To \QQ[a_1, a_2, x_2, x_3, x_5]$ which shows that $P/I$ is a free algebra.
\end{example}

Another case is obtained from Example 5.8.b. To be able to follow the steps of the algorithm better, 
we rename the setting of this example slightly.

\begin{example}\label{ex-similCrachiolaCont}
Let $P = \QQ[a, x_1, x_2, x_3]$ be graded by $W =(0,1,3,3)$,
consider $g = (1+a^2)\, x_2 + a x_3 + x_1^3$, and let $I=\langle g\rangle$.
Then~$I$ is a $W\!$-homogeneous ideal and $P/I$ is a regular ring.

Since the polynomial~$g$ is already a homogeneous regular sequence
and the entries of $(0\; 1{+}a^2\; a)$ generate the unit ideal in $\QQ[a]$,
the first two steps of Algorithm~\ref{alg-free} do nothing and we can immediately
start with step~(3). Let us follow the steps of Algorithm~\ref{alg-UMC-reemb}.

First of all, we have $D=\{1, 3\}$. In degree $d=1$, the matrix $A_1$ is empty 
and there is nothing to do. In degree $d=3$, we have $A_3 = (1{+}a^2\;\; a)$. 
The result of Algorithm~\ref{alg-completion} is the matrix
$B = \left( \begin{smallmatrix}  1 & -a\\  -a & a^2{+}1 \end{smallmatrix} \right)$.
This yields the $\QQ$-algebra automorphism $\phi:\; P \longrightarrow P$
defined by 
$$
F \;=\;  (a,\, x_1,\, -ax_3 +x_2,\, a^2x_3 -ax_2 +x_3)
$$
When we apply~$\phi$ to~$g$, we get $\phi(g) = x_2 - q$ with $q=-x_1^3$.

Thus we get $\widehat{F} = (a, x_1, -x_1^3, x_3)$ and the tuple
$T = (a, x_1, -x_1^3 - ax_3, a x_1^3 + (a^2+1) x_3)$ which defines a
$\QQ$-algebra isomorphism  $\Theta:\;  P/I \longrightarrow \QQ[a,x_1,x_3]$.
This shows that the ring $P/I$ is a free $\QQ$-algebra.
\end{example}

Our final two examples point out that the hypotheses of Algorithm~\ref{alg-UMC-reemb}
and Theorem~\ref{thm-Free} are indispensable.

\begin{example}\label{ex-fibernotequi}
Let $P =\QQ[a,x_1,x_2]$ be graded by $W=(0, 1, 1)$. We consider the ideal
$I = \langle (1-a)x_1, (1-2a)x_2 \rangle$. Then~$I$ is generated by a
homogeneous regular sequence, but the ring $P/I$ is not regular, as the
matrix $A= \left( \begin{smallmatrix} 1-a & 0 \\ 0 & 1-2a \end{smallmatrix} \right)$
fails to be unimodular. In particular, the $\QQ$-algebra $P/I$ is not free.
\end{example}

It is also important to note that the existence of a non-negative grading for which~$I$
is a homogeneous ideal is a crucial hypothesis. The next example demonstrates this
strikingly.

\begin{example}{\bf (The Circle)}\label{ex-circle}\\
Let $P = \QQ[x,y]$, and let $I = \langle x^2+y^2-2x \rangle$.
Clearly, the zero set $\mathcal{Z}(I)$ is a circle of radius~1 around $(1,0)$.
The ring $P/I$ is regular and a (global) complete intersection.
However, the ideal~$I$ is not homogeneous with respect to any non-trivial non-negative $\ZZ$-grading.
\begin{enumerate}
\item[(a)] There is no isomorphism between
$P/I$ and a univariate polynomial ring over~$\QQ$. To verify this, it suffices to note
that $\bar{x}$, $\bar{x}-1$, and $\bar{y}$ are irreducible, but not prime elements in~$P/I$. 
Therefore $P/I$ is not a factorial domain.

\item[(b)] It is also interesting to note that, if we enlarge the base field to $K=\QQ(i)$,
there is still no isomorphism between $P/I$ and $\QQ[z]$, but for a different reason.
In this case the $\QQ$-algebra isomorphism $\alpha:\; K[x,y] \longrightarrow K[u,v]$ defined by 
$\alpha(x) = {\tfrac{u+v}{2}} +1$, and $\alpha(y) = \tfrac{u-v}{2i}$ identifies~$I$ with
$\langle uv-1\rangle$. It induces an isomorphism $P/I \cong K[u,u^{-1}]$ with a factorial
domain, but $K[u,u^{-1}]$ is not isomorphic to a univariate polynomial ring over~$K$,
because it has too many invertible elements.
\end{enumerate}
\end{example}

\bigskip
\subsection*{Acknowledgements}
The second author thanks the University of Passau for its hospitality
during part of the preparation of this paper.

\bigskip\bigbreak
%
%


\begin{thebibliography}{99}


\bibitem{CoCoA} J.\ Abbott, A.M.\ Bigatti, and L.\ Robbiano,
CoCoA: a system for doing Computations in Commutative Algebra,
available at {\tt http://cocoa.dima.unige.it}.

\bibitem{BS} B.\ Barwick and B.\ Stone, Computing free bases for projective modules,
J.\ Softw.\ Algebra Geom.\ {\bf 5} (2013), 26--32.

\bibitem{BR} M.\ Beltrametti and L.\  Robbiano, Introduction to the
theory of weighted projective spaces, Expo. Math. {\bf 4} (1986), 111--162

\bibitem{Cr} A.J.\ Crachiola,
The hypersurface $x+x^2y+z^2+t^3$ over a field of arbitrary characteristic,
Proc. Amer. Math. Soc. {\bf 134} (2005), 1289--1298.


\bibitem{Eis} D.\ Eisenbud, {\it Commutative Algebra with a View Toward Algebraic Geometry},
Springer-Verlag, New York 1995


\bibitem{FQ} A.\ Fabianska and A.\ Quadrat, Application of the Quillen-Suslin theorem to
multidimensional systems theory, INRIA Research Report RR-6126, INRIA Rocquencourt, 
Le Chesnay-Rocquencourt, 2007.

\bibitem{F} N.\ Fitchas et al., Algorithmic aspects of Suslin's proof of Serre's
conjecture, Comput.\ Complexity {\bf 3} (1993), 31--55.

\bibitem{FG} N.\ Fitchas and A.\ Galligo, Nullstellensatz effectif et conjecture de Serre
(Theoreme de Quillen-Suslin) pour le calcul formel, Math.\ Nachr.\ {\bf 149} (1990), 
231--253.


\bibitem{Gup} N.\ Gupta, The Zariski cancellation problem and related problems 
in affine algebraic geometry, preprint 2022, 20 pages, 
{\tt https://doi.org/10.48550/arXiv.2208.14736}


\bibitem{KLR0} M.\ Kreuzer, L.N.~Long, and L.\ Robbiano, Computing subschemes of the
border basis scheme, Int. J. Algebra Comput. {\bf 30} (2020), 1671--1716.

\bibitem{KLR1} M.\ Kreuzer, L.N.~Long, and L.\ Robbiano, Cotangent spaces
and separating re-embeddings, J.\ Algebra Appl. {\bf 21} (2022), paper 2250188,
https://doi.org/10.1142/S0219498822501882

\bibitem{KLR2} M.\ Kreuzer, L.N.~Long, and L.\ Robbiano, Restricted Gr\"obner
fans and re-embeddings of affine algebras, S\~ao Paulo J.\ Math.\ Sci.\ {\bf 17} (2023), 
242--267,\\ {\tt https://doi.org/10.1007/s40863-022-00324-w}

\bibitem{KLR3} M.\ Kreuzer, L.N.~Long, and L.\ Robbiano, Re-embeddings of affine algebras
via Gr{\"o}bner fans of linear ideals, Beitr.\ Algebra Geom.\ (2024),
{\tt https://doi.org/10.1007/s13366-024-00733-2}

\bibitem{KR1} M.\ Kreuzer and L.\ Robbiano, {\it Computational
Commutative  Algebra 1}, Springer-Verlag, Berlin Heidelberg, 2000.

\bibitem{KR2} M.\ Kreuzer and L.\ Robbiano, {\it Computational
Commutative Algebra 2}, Springer-Verlag, Berlin Heidelberg, 2005.

\bibitem{KR3} M.\ Kreuzer and L.\ Robbiano, Deformations of border
bases, Collect.\ Math.\ {\bf 59} (2008), 275--297.

\bibitem{Kun} E.\ Kunz, {\it Introduction to Commutative Algebra and Algebraic Geometry},
Birkh\"auser, Boston 1985

\bibitem{LS} A.\ Logar and B.\ Sturmfels, Algorithms for the Quillen-Suslin theorem,
J.\ Algebra {\bf 145} (1992), 231--239.

\bibitem{RT} M.\ Roggero and L.\ Terracini, Ideals with an assigned initial ideal,
Int. Math. Forum \ {\bf 5} (2010), 2731--2750.

\bibitem{ApCoCoA} The ApCoCoA Team, ApCoCoA: Applied Computations in Computer Algebra,
available at {\tt apcocoa.uni-passau.de}


\end{thebibliography}
\end{document}